\documentclass[12pt, reqno]{amsart}
\usepackage{amssymb, mathrsfs}
\usepackage{latexsym}
\usepackage{amsmath}
\usepackage{amsthm}
\usepackage{amsfonts}
\usepackage{dsfont, upgreek}
\usepackage{mathtools}
\usepackage{epsfig}
\usepackage{amscd}
\usepackage{graphicx}
\usepackage{mathabx}

\usepackage[colorlinks=true,linkcolor=blue,citecolor=blue]{hyperref}

\usepackage{tikz-cd}

\usepackage{enumitem}
\usepackage{url}
\setlist[enumerate]{itemsep=0.5ex}

\usepackage{color}
\usepackage{tikz}
\usetikzlibrary{patterns}

\usepackage[left=1.2in,right=1.2in,top=1.4in,bottom=1.4in]{geometry}
\usetikzlibrary{decorations.markings}
\usetikzlibrary{arrows.meta}

\usepackage{extarrows}

\usepackage{adjustbox}
\usepackage{pgfplots}
\usepgfplotslibrary{colormaps}
\usepackage{slashed}

\usepackage[all,cmtip]{xy}
\usepackage{comment}

\theoremstyle{plain}
\newtheorem{theorem}{Theorem}[section]
\newtheorem{proposition}[theorem]{Proposition}
\newtheorem{lemma}[theorem]{Lemma}
\newtheorem{corollary}[theorem]{Corollary}
\newtheorem{question}[theorem]{Question}

\theoremstyle{definition} 
\newtheorem{definition}[theorem]{Definition}
\newtheorem{example}[theorem]{Example}
\newtheorem*{claim*}{Claim}

\theoremstyle{remark} 
\newtheorem{remark}[theorem]{Remark}

\numberwithin{equation}{section}

\newcommand{\N}{\mathbb{N}}
\newcommand{\diam}{\mathrm{diam}}
\newcommand{\Diam}{\mathrm{Diam}}
\newcommand{\m}{\mathrm{m}}

\newcommand{\supp}{\textnormal{supp }}
\newcommand{\In}{\textnormal{In}}

\newcommand{\ad}{\textnormal{ad}}
\newcommand{\dad}{\textnormal{dad}}
\newcommand{\asdim}{\textnormal{asdim}}
\newcommand{\dasdim}{\textnormal{dasdim}}

\newcommand{\interior}[1]{%
	{\kern0pt#1}^{\mathrm{\,o}}%
}

\makeatletter
\let\save@mathaccent\mathaccent
\newcommand*\if@single[3]{%
	\setbox0\hbox{${\mathaccent"0362{#1}}^H$}%
	\setbox2\hbox{${\mathaccent"0362{\kern0pt#1}}^H$}%
	\ifdim\ht0=\ht2 #3\else #2\fi
}
\newcommand*\rel@kern[1]{\kern#1\dimexpr\macc@kerna}
\newcommand*\wideaccent[2]{\@ifnextchar^{{\wide@accent{#1}{#2}{0}}}{\wide@accent{#1}{#2}{1}}}
\newcommand*\wide@accent[3]{\if@single{#2}{\wide@accent@{#1}{#2}{#3}{1}}{\wide@accent@{#1}{#2}{#3}{2}}}
\newcommand*\wide@accent@[4]{%
	\begingroup
	\def\mathaccent##1##2{%
		\let\mathaccent\save@mathaccent
		\if#42 \let\macc@nucleus\first@char \fi
		\setbox\z@\hbox{$\macc@style{\macc@nucleus}_{}$}%
		\setbox\tw@\hbox{$\macc@style{\macc@nucleus}{}_{}$}%
		\dimen@\wd\tw@
		\advance\dimen@-\wd\z@
		\divide\dimen@ 3
		\@tempdima\wd\tw@
		\advance\@tempdima-\scriptspace
		\divide\@tempdima 10
		\advance\dimen@-\@tempdima
		\ifdim\dimen@>\z@ \dimen@0pt\fi
		\rel@kern{0.6}\kern-\dimen@
		\if#41
		#1{\rel@kern{-0.6}\kern\dimen@\macc@nucleus\rel@kern{0.4}\kern\dimen@}%
		\advance\dimen@0.4\dimexpr\macc@kerna
		\let\final@kern#3%
		\ifdim\dimen@<\z@ \let\final@kern1\fi
		\if\final@kern1 \kern-\dimen@\fi
		\else
		#1{\rel@kern{-0.6}\kern\dimen@#2}%
		\fi
	}%
	\macc@depth\@ne
	\let\math@bgroup\@empty \let\math@egroup\macc@set@skewchar
	\mathsurround\z@ \frozen@everymath{\mathgroup\macc@group\relax}%
	\macc@set@skewchar\relax
	\let\mathaccentV\macc@nested@a
	\if#41
	\macc@nested@a\relax111{#2}%
	\else
	\def\gobble@till@marker##1\endmarker{}%
	\futurelet\first@char\gobble@till@marker#2\endmarker
	\ifcat\noexpand\first@char A\else
	\def\first@char{}%
	\fi
	\macc@nested@a\relax111{\first@char}%
	\fi
	\endgroup
}
\makeatother

\makeatletter
\newcommand*{\transpose}{%
	{\mathpalette\@transpose{}}%
}
\newcommand*{\@transpose}[2]{%
	\raisebox{\depth}{$\m@th#1\intercal$}%
}
\makeatother

\begin{document}
	
	\title[\tiny Dynamic asymptotic dimension growth for group actions and groupoids]{Dynamic asymptotic dimension growth for group actions and groupoids}

\author{Hang Wang}
\address[Hang Wang]{\normalfont{Research Center for Operator Algebras, School of Mathematical Sciences, East China Normal University, Shanghai 200062, China}}
\email{wanghang@math.ecnu.edu.cn}
\thanks{Hang Wang is supported by NSFC (Grant No. 12271165) and in part by the Science and Technology Commission of Shanghai Municipality (Grant Nos. 22DZ2229014, 23JC1401900). }

\author{Yanru Wang}
\address[Yanru Wang]{\normalfont{Tianyuan Mathematical Center in Southwest China, School of Mathematics\\ Sichuan University, Chengdu 610065, China}}
\email{yanruwang@scu.edu.cn}
\thanks{Yanru Wang is supported by the Science and Technology Commission of Shanghai Municipality (Grant No. 23JC1401900). Yanru Wang is the corresponding author.}

\author{Jianguo Zhang}
\address[Jianguo Zhang]{\normalfont{School of Mathematics and Statistics, Shaanxi Normal University,
Xi'an 710119, China}}
\email{jgzhang@snnu.edu.cn}
\thanks{Jianguo Zhang is supported by NSFC (Grant Nos. 12271165, 12171156, and 12301154). }

\author{Dapeng Zhou}
\address[Dapeng Zhou]{\normalfont{School of Statistics and Information, Shanghai University of International Business and Economics, Shanghai 201620, China}}
\email{dpzhou@suibe.edu.cn}
\thanks{Dapeng Zhou is supported by NSFC (Grant No. 12271165).}

	
 \begin{abstract}
We introduce the notion of dynamic asymptotic dimension growth for actions of discrete groups on compact spaces, and more generally for locally compact \'etale groupoids. Moreover, we demonstrate that the asymptotic dimension growth for a discrete metric space of bounded geometry is equivalent to the dynamic asymptotic dimension growth for its associated coarse groupoid. Consequently, we deduce that the coarse groupoid with subexponential dynamic asymptotic dimension growth is amenable. More generally, we show that every $\sigma$-compact locally compact Hausdorff \'etale groupoid with compact unit space and dynamic asymptotic dimension growth at most $x^{\alpha}$ $(0<\alpha<1)$ is amenable. As an application, we show that the Baum--Connes conjecture with coefficients holds for such groupoids.
 \end{abstract}

	\maketitle
\tableofcontents

\section{Introduction}
The notion of asymptotic dimension was primarily proposed by Gromov in \cite{Gro93} as a large-scale analog of the Lebesgue covering dimension. This concept gained widespread attention following Yu's seminal work, which established that the Novikov conjecture holds for groups with finite asymptotic dimension (see \cite{Yu98}). This pivotal result spurred extensive research into groups and spaces with finite asymptotic dimension. For instance, Gromov initially established that hyperbolic groups and metric spaces have finite asymptotic dimension (see \cite{Gro93}). A more straightforward proof of this result was later provided by Roe (see \cite{Roe05}). Building on these foundations, Bell and Dranishnikov showed that the free product, the amalgamated free product, and the HNN extension can preserve finite asymptotic dimension (see \cite{BD01}). Four years later, Osin demonstrated that if a finitely generated group $G$ is hyperbolic relative to a collection of subgroups with finite asymptotic dimension, then $G$ itself has finite asymptotic dimension (see \cite{Osi05}). Further advancements were made by Wright, who proved that the asymptotic dimension of a finite-dimensional CAT(0) cube complex is bounded above by its dimension (see \cite{Wri12}). Additionally, Bestvina, Bromberg, and Fujiwara established that mapping class groups have finite asymptotic dimension (see \cite{BBF15}). Following this, Behrstock, Hagen, and Sisto proved that all hierarchically hyperbolic groups have finite asymptotic dimension---a result that subsumes earlier findings for mapping class groups, relatively hyperbolic groups, and CAT(0) cube complexes (see \cite{BHS17}), among many others. For a comprehensive survey of this theory, refer to \cite{BD11}; for recent developments, see \cite{FP21, Tse23, FS24}. 

Moreover, finite asymptotic dimension has found applications in the domain of $L^{p}$ operator algebras as well. Analogous to the coarse Baum--Connes conjecture, the third and fourth authors established the $L^{p}$ coarse Baum--Connes conjecture for $p\in [1,\infty)$ and proved that this conjecture holds for proper metric spaces with finite asymptotic dimension (see \cite{ZZ21}). Subsequent work has focused on discrete metric spaces, where the  hypothesis of bounded geometry becomes essential. Recall that a discrete metric space $(X,d)$ has bounded geometry if for every $r>0$, there exists a constant $n_{r}\in\mathbb{N}$ such that every ball $B(x,r)$ has cardinality at most $n_{r}$. Within this framework, the present authors first investigated the persistence approximation property in the quantitative $K$-theory of filtered $L^{p}$ operator algebras and demonstrated that any $L^p$ Roe algebra for a coarsely uniformly contractible discrete metric space with bounded geometry and finite asymptotic dimension satisfies this property (see \cite{WWZZ24}). Subsequently, the authors provided an alternative proof of the $L^{p}$ coarse Baum--Connes conjecture for discrete metric spaces with bounded geometry and finite asymptotic dimension via $C_{0}$ coarse geometry (see \cite{WWZZ25}). 

However, numerous groups are known to have infinite asymptotic dimension, including the wreath product $\mathbb{Z}\wr\mathbb{Z}$ (see \cite{Dra06}), the first Grigorchuk group (see \cite{Smi07}), Thompson's group (see \cite{DS11}), and Gromov's celebrated group in \cite{Gro03} containing an expander. To study the growth type of dimension functions of groups with infinite asymptotic dimension, Dranishnikov \cite{Dra06} introduced the notion of asymptotic dimension growth for metric spaces. This concept extends the notion of finite asymptotic dimension, since a constant asymptotic dimension function indicates finite asymptotic dimension. Dranishnikov also established that the wreath product of a finitely generated nilpotent group with a finitely generated group of finite asymptotic dimension exhibits polynomial asymptotic dimension growth. Furthermore, he proved that polynomial asymptotic dimension growth implies Yu's property A \cite{Dra06}, a weak amenability-type condition introduced by Yu \cite{Yu00}. In conjunction with Yu's theorem on the coarse Baum--Connes conjecture \cite{Yu00}, it follows that spaces with polynomial asymptotic dimension growth satisfy the coarse Baum--Connes conjecture, and consequently, the Novikov conjecture. Building on Dranishnikov's insights, Ozawa showed that a metric space with subexponential asymptotic dimension growth has property A \cite{Oza12}; see also \cite{Opp14} for the case of discrete metric spaces of bounded geometry. In this context, subexponential asymptotic dimension growth is a property that is weaker than finite asymptotic dimension but stronger than property A. 

Many examples of groups and spaces exhibit subexponential asymptotic dimension growth. A fascinating example is a geodesic coarse median space with finite rank and at most exponential volume growth, which has subexponential asymptotic dimension growth \cite{ANWZ18}. This space includes the hyperbolic space and the mapping class group. Another intriguing example is the group $G=\bigoplus\limits_{n\in\mathbb{N}}f(n)\mathbb{Z}$, where  $f:\mathbb{N}\rightarrow\mathbb{R}^{+}$ is a strictly increasing function that grows faster than a linear function. When endowed with the direct sum metric, $G$ becomes a discrete metric space of bounded geometry with slow asymptotic dimension growth $f^{-1}$ \cite{CFY08}---a special case of subexponential asymptotic dimension growth. Additionally, Bell examined how various constructions involving groups affect the asymptotic dimension function \cite{Bel05}. In this paper, we introduce four equivalent interlocking definitions of asymptotic dimension growth. The second definition, Definition \ref{def 2.4}, is derived from [\cite{ANWZ18}, Lemma 2.6]; the fourth, Definition \ref{def 2.15}, is formulated in the language of coarse geometry and is applied in the proof of Theorem \ref{th 1.1}.

Inspired by Gromov’s theory of asymptotic dimension, Guentner, Willett, and Yu introduced the dynamic asymptotic dimension in their seminal paper \cite{GWY17}, which is a concept of dimension for actions of discrete groups on locally compact spaces, as well as locally compact Hausdorff \'etale groupoids. Since its inception, this notion has found various applications. It provides upper bounds on the nuclear dimension of the resulting groupoid $C^{*}$-algebras \cite{GWY17, DS18, ALSS21, HW24} and twisted groupoid $C^{*}$-algebras \cite{CDGHV24, BL24}. Additionally, it establishes a sufficient condition for Matui's HK conjecture in low dimensions \cite{BDGW23}. Beyond its applications, the dynamic asymptotic dimension is intimately connected to other key concepts in the field. It correlates closely with the diagonal dimension of sub-$C^{*}$-algebras \cite{LLW23}, Kerr's tower dimension \cite{Ker20}, and the conditions formulated by Bartels, Lück, and Reich \cite{BLR08, Bar16}. In a recent study \cite{Bon24}, B\"onicke furthered the study of this dimension by establishing several permanence properties including estimates for products and unions of groupoids,  as well as comparing the asymptotic dimension of the resulting coarse space with the dynamic asymptotic dimension of the underlying groupoid. 

While finite dynamic asymptotic dimension has found numerous applications, the study of group actions and groupoids with infinite dynamic asymptotic dimension remains largely unexplored. The aim of this paper is to contribute to their study by examining the growth of the dynamic asymptotic dimension function for these group actions and groupoids. Our work draws inspiration from three key references: Guentner, Willett, and Yu's foundational research on finite dynamic asymptotic dimension for group actions and \'etale groupoids in \cite{GWY17}; their subsequent investigation verifying the Baum--Connes conjecture for actions with finite dynamical complexity in \cite{GWY24}; and the contributions of Arzhantseva, Niblo, Wright, and Zhang on characterizing asymptotic dimension growth in \cite{ANWZ18}. These influential works collectively provide valuable methods and results for our investigation into the infinite-dimensional case.

In the case of group actions, we introduce the concept of dynamic asymptotic dimension growth for actions of discrete groups on compact Hausdorff spaces in Definition \ref{def 3.2}. This notion serves as a dynamical counterpart of asymptotic dimension growth, and as a generalization of finite dynamic asymptotic dimension. Building upon existing estimates of the dynamic asymptotic dimension for various group actions \cite{ALSS21, CJMST23, GWY17}, we subsequently construct examples of group actions with exponential dynamic asymptotic dimension growth and infinite dynamic asymptotic dimension. This is achieved by applying Lemma \ref{lemma 4.13} and Theorem \ref{th 1.1} to examples of groups with exponential asymptotic dimension growth.

In the case of \'etale groupoids, we define the concept of dynamic asymptotic dimension growth for locally compact Hausdorff \'etale groupoids, thereby extending our earlier notion for actions of discrete groups on compact spaces. Specifically, we first define the dynamic asymptotic dimension growth of an \'etale groupoid with compact unit space in Definition \ref{def 4.6}. For a groupoid with non-compact unit space, our definition reduces to the compact case via the Alexandrov groupoid defined in [\cite{KKLRU21}, Definition 7.8].  Within the framework of coarse geometry, we establish a connection between our concept and the asymptotic dimension growth, leading to the following theorem:

\begin{theorem}[see Theorem \ref{th 4.17}]\label{th 1.1}
Let $X$ be a discrete metric space of bounded geometry, and let $G(X)$ be the associated coarse groupoid. Let $f: \mathbb{R}^{+}\rightarrow\mathbb{N}$ be a non-decreasing function. Then $\ad_{X}\approx f$ if and only if $\dad_{G(X)}\approx f$. 
\end{theorem}

This theorem generalizes [\cite{GWY17}, Theorem 6.4]. It serves as our principal motivation due to its connection to the coarse Baum--Connes conjecture. The pivotal step in our proof involves a combination of the fourth equivalent definition of asymptotic dimension growth with the technical approach underlying [\cite{GWY17}, Theorem 6.4]. Furthermore, this theorem provides a shortcut to construct concrete examples of group actions with at most exponential dynamic asymptotic dimension growth and infinite dynamic asymptotic dimension, as exemplified in Example \ref{example 3.6} and  Example \ref{example 3.7}.

These insights lead us to the central theme of the paper: the amenability of groupoids. Amenability plays a fundamental role in the field of operator algebras and has led to significant applications. A key result by Tu establishes that the Baum--Connes conjecture with coefficients holds for amenable groupoids \cite{Tu99}. Consequently, the reduced $C^{*}$-algebras of amenable groupoids satisfy the universal coefficient theorem (UCT) of Rosenberg and Schochet \cite{Tu99}. This further implies that for $C^{*}$-algebras arising from amenable groupoids, if they are separable, simple, and $\mathcal{Z}$-stable, then they are classified by their Elliott invariants \cite{EGLN25, GL24, Li20, GLN20, KP00, Phi00, TWW17}. Furthermore, the rational HK conjecture is satisfied by second-countable, amenable, principal, ample, and \'etale groupoids \cite{PY26}. These results motivate the following central question for the remainder of the paper:

\begin{question}\label{q 1.2}
What kind of dynamic asymptotic dimension growth of \'etale groupoids implies amenability? 
\end{question}

Returning to the context of coarse groupoids, which are a class of \'etale groupoids, we obtain the following corollary providing an alternative sufficient condition for amenability. 

\begin{corollary}[see Corollary \ref{cor 4.21}]\label{cor 1.3}
Let $X$ be a discrete metric space of bounded geometry. If the associated coarse groupoid $G(X)$ has subexponential dynamic asymptotic dimension growth, then $G(X)$ is amenable.
\end{corollary}

This result follows from the established equivalence between property A for a discrete metric space $X$ of bounded geometry and the topological amenability of the coarse groupoid $G(X)$ [\cite{STY02}, Theorem 5.3]. By combining this with Theorem \ref{th 1.1} and [\cite{Oza12}, Theorem 1], we reformulate [\cite{Oza12}, Theorem 1] in the language of groupoids. We now illustrate this result with a concrete example: for a geodesic coarse median space $X$ of bounded geometry with finite rank and at most exponential volume growth, the corresponding coarse groupoid $G(X)$ has subexponential dynamic asymptotic dimension growth, thereby confirming its amenability.

For a general \'etale groupoid $\mathcal{G}$, the following main theorem establishes a sufficient condition for its amenability, thereby providing a positive answer to Question \ref{q 1.2}.

\begin{theorem}(see Theorem \ref{th 6.4})\label{th 1.4}
Let $\mathcal{G}$ be a $\sigma$-compact, locally compact,  Hausdorff, and \'etale groupoid with compact unit space and dynamic asymptotic dimension growth $f\preceq x^{\alpha} (0<\alpha<1)$. Then $\mathcal{G}$ is amenable.
\end{theorem}

The proof employs generalized partitions of unity (Proposition \ref{prop 5.1}) together with the approach developed by Guentner, Willett, and Yu for studying finite dynamical complexity in \'etale groupoids [\cite{GWY24}, Theorem A.9]. Our main theorem extends Corollary \ref{cor 1.3} to general \'etale groupoids with slow growth. Using a similar approach, we also show that every locally compact Hausdorff \'etale groupoid with finite dynamic asymptotic dimension is amenable (Theorem \ref{th 6.3}), thereby strengthening the earlier result in [\cite{GWY17}, Corollary 8.25] by removing the freeness assumption.

By combining Theorem \ref{th 1.4} with the result from \cite{Tu99}, we conclude that the Baum--Connes conjecture with coefficients holds for $\sigma$-compact, locally compact, Hausdorff, and \'etale groupoids with compact unit space and dynamic asymptotic dimension growth at most $x^{\alpha}$ $(0<\alpha<1)$. Furthermore, applying this conclusion to  [\cite{PY26}, Corollary 5.6], we deduce that the rational HK conjecture holds for second-countable, $\sigma$-compact, locally compact, Hausdorff, principal, ample, and \'etale groupoids with compact unit space and dynamic asymptotic dimension growth at most $x^{\alpha}$ $(0<\alpha<1)$.

The paper is organized as follows. In Section 2, we review the concept of asymptotic dimension growth and present four equivalent definitions of this concept. In Section 3, we introduce the notion of dynamic asymptotic dimension growth for group actions on compact Hausdorff spaces. In Section 4, we extend the definition of dynamic asymptotic dimension growth to \'etale groupoids and provide several concrete examples. Furthermore, we demonstrate that the asymptotic dimension growth for a discrete metric space of bounded geometry is equivalent to the dynamic asymptotic dimension growth for its associated coarse groupoid. In Section 5, we introduce a generalization of partitions of unity, which serves as a key technical tool for establishing the amenability of groupoids. In Section 6, we establish a sufficient condition for the amenability of groupoids: every $\sigma$-compact, locally compact, Hausdorff, and \'etale groupoid $\mathcal{G}$ with compact unit space and dynamic asymptotic dimension growth at most $x^{\alpha}$ $(0<\alpha<1)$ is amenable. This result leads to applications involving the Baum--Connes conjecture and Matui's HK conjecture.

\section{Asymptotic dimension growth for metric spaces}\label{section 2}
The volume growth of finitely generated groups is at most exponential. To investigate the intermediate dimension growth of these groups, Dranishnikov proposed the concept of asymptotic dimension growth in \cite{Dra06}, which generalizes the notion of finite asymptotic dimension. In this section, we present four equivalent definitions of asymptotic dimension growth. The fourth equivalent definition, formulated in the language of coarse geometry, plays a crucial role in the proof of Theorem \ref{th 4.17} in Section \ref{section: dad}. To illustrate a more intuitive understanding of this concept, we provide examples of spaces exhibiting exponential growth, subexponential growth, polynomial growth, and slow growth.
\subsection{Equivalent definitions}
In this subsection, we review the concept of asymptotic dimension growth and present three equivalent definitions. To enhance clarity, we illustrate this concept with several concrete examples. 

To clarify the definition of the asymptotic dimension function, we outline some elementary concepts. Let $(X,d)$ be a metric space and $R>0$. A family $\mathcal{W}$ of subsets of $X$ is termed $R$-disjoint if the distance $d(W, W')>R$ for every pair of distinct subsets $W\neq W'$ within $\mathcal{W}$. A family $\mathcal{V}$ is termed uniformly bounded if its diameter, defined as $\Diam(\mathcal{V})=\sup\{\diam(V)\mid V\in\mathcal{V}\}$, is finite.

Let $\mathcal{U}$ be a cover of $X$. The Lebesgue number of $\mathcal{U}$, denoted by $L(\mathcal{U})$, is the maximum number $R$ such that for every subset $V$ of $X$ with diameter $\diam(V)\leq R$, there exists a set $U\in\mathcal{U}$ such that $V\subseteq U$. The multiplicity of $\mathcal{U}$, denoted as $\m(\mathcal{U})$, is the largest number of elements of $\mathcal{U}$ with nonempty intersection. The $R$-multiplicity of $\mathcal{U}$, denoted by $\m_{R}(\mathcal{U})$, is the smallest integer $n$ such that for any $x\in X$, the ball $B(x,R)$ intersects at most $n$ sets of $\mathcal{U}$. 

We now proceed to recall the definition of the asymptotic dimension function, originally given by Dranishnikov in \cite{Dra06}.
\begin{definition}[\cite{Dra06}]
Let $(X,d)$ be a metric space. For any $R>0$, let $\ad_{X}(R)$ be the smallest number $n\in\mathbb{N}$ for which there exists a uniformly bounded cover $\mathcal{U}=\{U_{i}\}_{i\in I}$ of $X$ with $L(\mathcal{U})>R$ and $\m(\mathcal{U})\leq n+1$. We call $\ad_{X}: \mathbb{R}^{+}\rightarrow\mathbb{N}$ the asymptotic dimension function of $X$.
\end{definition}
\begin{remark}
The following statements hold:
\begin{enumerate}
\item The function $\ad_{X}$ is nondecreasing, and the asymptotic dimension of $X$ is defined as follows: 
\begin{equation*}
\asdim(X)=\lim\limits_{R\rightarrow\infty}\ad_{X}(R).
\end{equation*}
\item $\asdim (X)\leq d$ if and only if $\ad_{X}$ is constant and $\ad_{X}\leq d$.
\end{enumerate}
\end{remark}

Analogous to the growth type of the volume function, the growth type of the asymptotic dimension function is an essential concept. 

Recall from \cite{Bel05} that for $f,g:\mathbb{R}^{+}\rightarrow\mathbb{R}^{+}$,
\begin{itemize}
\item we say that $f$ is at most $g$, denoted as $f\preceq g$, if there exists $k\in\mathbb{N}$ such that $f(x)\leq kg(kx+k)+k$ for any $x>k$;
\item we say that $f$ is equivalent to $g$, denoted as $f\approx g$, if both $f\preceq g$ and $g\preceq f$. Clearly, $\approx$ is an equivalence relation;
\item we define the growth type of $f$ to be the $\approx$-equivalence class of $f$;
\item we define the \textit{asymptotic dimension growth} of $X$ to be the growth type of $\ad_{X}$.
\end{itemize}

The following proposition establishes a key property of the asymptotic dimension growth.
Although the asymptotic dimension growth is not a coarse invariant like the finite asymptotic dimension, it turns out to be a quasi-isometric invariant by the following result \cite{Bel05, Dra06}. Thus, it is a group invariant because any two word metrics on a finitely generated group are quasi-isometric.

\begin{proposition}[\cite{Bel05, Dra06}]
Let $X$ and $Y$ be two discrete metric spaces with bounded geometry. If $X$ and $Y$ are quasi-isometric, then $\ad_{X}\approx \ad_{Y}$. In particular, the asymptotic dimension growth is well-defined for finitely generated groups.
\end{proposition}

We now present the second equivalent definition that characterizes the asymptotic dimension growth by the existence of a cover satisfying the $R$-multiplicity condition. Crucially, the following definition is shown to be equivalent to the original one, as established in the subsequent lemma.

\begin{definition}[\cite{ANWZ18}]\label{def 2.4}
Let $(X,d)$ be a metric space. For any $R>0$, we define $\widetilde{\ad}_{X}(R)$ to be the smallest number $n\in\mathbb{N}$ for which there exists a uniformly bounded cover $\mathcal{V}=\{V_{j}\}_{j\in J}$ with $\m_{R}(\mathcal{V})\leq n+1$.
\end{definition}

\begin{lemma}[\cite{ANWZ18}]\label{lemma 2.5}
Let $(X,d)$ be a metric space. Then $\widetilde{\ad}_{X}\approx\ad_{X}$.
\end{lemma}

The third equivalent definition reframes the asymptotic dimension growth in the context of a separable metric space $X$ through a cover consisting of $R$-disjoint families.

\begin{definition}\label{def 2.6}
Let $(X,d)$ be a separable metric space. For any $R>0$, we define $\widehat{\ad}_{X}(R)$ to be the smallest number $n\in\mathbb{N}$ for which there exist uniformly bounded, $R$-disjoint families $\mathcal{W}_{0},\dots,\mathcal{W}_{n}$ of subsets of $X$ such that $\bigcup\limits_{l=0}^{n}\mathcal{W}_{l}$ is a cover of $X$.
\end{definition}

The following proposition shows that this new formulation is equivalent to the previous two definitions. This equivalence is a pivotal step, as it directly enables the transition to the fourth definition formulated in the language of coarse geometry in the next subsection.

\begin{proposition}\label{prop 2.7}
Let $(X,d)$ be a separable metric space. Then $\widehat{\ad}_{X}\approx\widetilde{\ad}_{X}\approx\ad_{X}$.
\end{proposition}
\begin{proof}
Given $R>0$, it follows that $2R>R$. Suppose there exist uniformly bounded, $2R$-disjoint families $\mathcal{W}_{0},\dots,\mathcal{W}_{\widehat{\ad}_{X}(2R)}$ of subsets of $X$ that cover $X$. Put $\mathcal{W}=\bigcup\limits_{l}\mathcal{W}_{l}$, which is obviously a uniformly bounded cover. We claim that $\m_{R}(\mathcal{W})\leq\widehat{\ad}_{X}(2R)$. By the definition of $R$-multiplicity, for any $x\in X$, there are at most $\widehat{\ad}_{X}(2R)+1$ sets in $\mathcal{W}$ having non-empty intersection with the closed ball $B(x,R)$. Otherwise, there must be two sets $W$ and $W'$ from the same family $\mathcal{W}_{i}$ for some $i$ such that both $W\cap B(x, R)\neq\varnothing$ and $W'\cap B(x,R)\neq\varnothing$, since there are only $\widehat{\ad}_{X}(2R)+1$ different families. This implies that $d(W,W')\leq 2R$, contradicting the $2R$-disjoint property of the family $\mathcal{W}_{i}$. Therefore, $\m_{R}(\mathcal{W})\leq\widehat{\ad}_{X}(2R)$ as claimed, and thus yields $\widetilde{\ad}_{X}\preceq\widehat{\ad}_{X}$.

Conversely, given $R>0$, suppose $\mathcal{V}=\{V_{j}\}_{j\in J}$ is a unifromly bounded cover with $\m_{R}(\mathcal{V})\leq \widetilde{\ad}_{X}(R)+1$. Since $X$ is separable, it is second-countable, which implies that any cover of $X$ admits a countable subcover. Therefore, we assume that the cover $\mathcal{V}$ is countable. We initially claim that there exist uniformly bounded, $R$-disjoint families $\mathcal{V}_{0},\dots,\mathcal{V}_{\widetilde{\ad}_{X}(R)}$ of subsets of $X$ that cover $X$. This claim is demonstrated through a three-step process.

\textbf{Step 1}. In the first step, we construct an irreducible cover of $X$. Let $\mathcal{U}=\{U_{i}\}_{i\in I}$ be a subcover of $\mathcal{V}$ such that 
$$
V_{j}\bigg\backslash\left(\bigcup\limits_{V'\in\mathcal{V}\backslash\{V_{j}\}}V'\right)\neq\varnothing, \text{ for all } V_{j}\in\mathcal{V},
$$
which forms an irreducible cover of $X$. As $\m_{R}(\mathcal{U})\leq \m_{R}(\mathcal{V})\leq\widetilde{\ad}_{X}(R)+1$, we have $\m_{R}(\mathcal{U})\leq\widetilde{\ad}_{X}(R)+1$ ( where the equality can be reached here ).

\textbf{Step 2}. In this step, we build $\widetilde{\ad}_{X}(R)+1$ boxes. Since $\m_{R}(\mathcal{U})\leq \widetilde{\ad}_{X}(R)+1$ and the equality is achievable, there exists $x\in X$ such that the ball $B(x,R)$ intersects $\widetilde{\ad}_{X}(R)+1$ elements $U_{0},\dots, U_{\widetilde{\ad}_{X}(R)}$ of $\mathcal{V}$. Consequently, $d(U_{i}, U_{j})\leq 2R$, for all $i,j\in\{0,\dots,\widetilde{\ad}_{X}(R)\}$ with $i\neq j$. We then put $U_{0},\dots,U_{\widetilde{\ad}_{X}(R)}$ into $\widetilde{\ad}_{X}(R)+1$ distinct boxes. 

\textbf{Step 3}. In the final step, we allocate the remaining sets of the cover $\mathcal{U}$ into these $\widetilde{\ad}_{X}(R)+1$ boxes. Set $\mathcal{U}\big\backslash\bigcup\limits_{i=0}^{\widetilde{\ad}_{X}(R)}U_{i}=\{U_{k}\}_{k\in I}$, and let $\Box_{l}$ represent the $l$-th box. We next claim that for any $x\in U_{k}$, there exists $l\in\{0,\dots,\widetilde{\ad}_{X}(R)\}$ such that 
$$
d(x,B_{l})>R, \forall B_{l}\in \Box_{l}.
$$
Otherwise, there exists $y\in U_{k}$ such that for all $l\in\{0,\dots,\widetilde{\ad}_{X}(R)\}$, there is a set $B_{l}\in\Box_{l}$ such that 
$$
d(y, B_{l})\leq R,
$$
which implies that the ball $B(y, R)$ intersects $\widetilde{\ad}_{X}(R)+2$ sets $U_{k},B_{0},\dots, B_{\widetilde{\ad}_{X}(R)}$ of $\mathcal{U}$, which  contradicts that $\m_{R}(\mathcal{U})\leq \widetilde{\ad}_{X}(R)+1$. Therefore, for any $x\in U_{k}$, $x$ must belong to some box. Subsequently, we divide $U_{k}$ into several pieces and place them into boxes sequentially. Let 
\begin{align*}
 U^{(0)}_{k}&=\{x\in U_{k}\mid d(x, B_{0})>R, \text{ for all }B_{0}\in\Box_{0}\},\\
U^{(l)}_{k}&=\{x\in U_{k}\big\backslash\bigcup\limits_{n=0}^{l-1}U^{(n)}_{k}\mid d(x, B_{l})>R, \text{ for all }B_{l}\in\Box_{l}\}, l=1,\dots,\widetilde{\ad}_{X}(R).
\end{align*}
Then we put $U^{(l)}_{k}$ into the box $\Box_{l}$ for $l=0,\dots,\widetilde{\ad}_{X}(R)$ in ascending order. If $U^{(l)}_{k}$ is the empty set, we proceed to the next box. Observe that $U^{(l)}_{k}$ is a subset of $U_{k}$, and $U_{k}=\bigsqcup\limits_{l=0}^{\widetilde{\ad}_{X}(R)}U^{(l)}_{k}$. Similarly, we can divide all elements of $\mathcal{U}\big\backslash\bigcup\limits_{j=0}^{\widetilde{\ad}_{X}(R)}U_{j}$ into $\widetilde{\ad}_{X}(R)+1$ boxes. Consequently, we obtain new boxes and for all $l\in\{0,\dots,\widetilde{\ad}_{X}(R)\}$, 
$$
d(B_{l},B'_{l})>R, \text{ for every }B_{l}\neq B'_{l} \text{ in }\Box_{l}.
$$
Let $\mathcal{W}_{l}=\{B_{l}\mid B_{l}\in\Box_{l}\}$, then each $\mathcal{W}_{l}$ is a uniformly bounded, $R$-disjoint family and  $X=\bigcup\limits_{l=0}^{\widetilde{\ad}_{X}(R)}\mathcal{W}_{l}$. By Definition \ref{def 2.6}, we have $\widehat{\ad}_{X}(R)\leq \widetilde{\ad}_{X}(R)$, which then yields $ \widehat{\ad}_{X}\preceq \widetilde{\ad}_{X}$. Hence, $\widehat{\ad}_{X}\approx\widetilde{\ad}_{X}$.

By the transitivity of the equivalence relation, we conclude that $\widehat{\ad}_{X}\approx\widetilde{\ad}_{X}\approx\ad_{X}$ by Lemma \ref{lemma 2.5}.
\end{proof}

In the second part of this subsection, we present several typical examples of groups and spaces exhibiting asymptotic dimension growth. We primarily focus on four growth types: exponential, subexponential, polynomial, and slow. We define a monotone increasing function $f:\mathbb{R}^{+}\rightarrow\mathbb{R}^{+}$ as having subexponential growth if for any $b>1$,
$$
\lim\limits_{t\rightarrow\infty}\frac{f(t)}{b^{t}}=0.
$$
Observe that if $f$ has polynomial growth, it inherently exhibits subexponential growth, since polynomial functions grow at a slower rate than exponential functions. A space $X$ is said to have polynomial (resp. subexponential, exponential) asymptotic dimension growth, if there exists a polynomial (resp. subexponential, exponential) function $f$ such that $\ad_{X}\preceq f$. Recall from \cite{Dra00} that $X$ has slow asymptotic dimension growth, if there exists a linear function $f$ such that $\ad_{X}\preceq f$.

Similar to volume growth, a finitely generated group has at most exponential asymptotic dimension growth.
\begin{example}[\cite{Dra06}]\label{example 2.8}
Let $\Gamma$ be a finitely generated group equipped with the word metric. Under this metric, $\Gamma$ becomes a discrete metric space since the distance between any two distinct elements is at least one.  
Furthermore, $\Gamma$ has bounded geometry since the finite generating set induces a Cayley graph that is uniformly locally finite. Moreover, there exists $\alpha>0$ such that $\ad_{\Gamma}(R)\preceq e^{\alpha R}$. 
\end{example}

Ozawa has demonstrated that subexponential asymptotic dimension growth implies Yu's property A \cite{Oza12}, thereby identifying two classes of spaces having this property as follows. 

\begin{example}[\cite{ANWZ18}]\label{example 2.9}
Let $X$ be a geodesic coarse median space with finite rank and at most exponential volume growth, then $X$ has subexponential asymptotic dimension growth. 
\end{example}

\begin{example}[\cite{Dra06}]\label{example 2.10}
Let $N$ be a finitely generated nilpotent group and $G$ be a finitely generated group with finite asymptotic dimension. Then the wreath product $N\wr G$ has polynomial asymptotic dimension growth. In particular, $\mathbb{Z}\wr\mathbb{Z}$ has polynomial asymptotic dimension growth.
\end{example}

The subsequent example demonstrates that given a strictly increasing function $f$, we can create a metric space with asymptotic dimension growth $f^{-1}$. This method provides an approach to construct a metric space with slow asymptotic dimension growth by choosing a function $f$ that grows faster than a linear function.

\begin{example}[\cite{CFY08}]\label{example 2.11}
Let $f:\mathbb{N}\rightarrow\mathbb{R}^{+}$ be a strictly increasing function. Consider the direct sum
$G=\bigoplus\limits_{n\in\mathbb{N}}f(n)\mathbb{Z}$, and endow it with the direct sum metric defined by
\begin{equation}\label{eq 2.2}
d(g,h)=\sup\limits_{n\in\mathbb{N}}|g_{n}-h_{n}|,
\end{equation}
where $g=(g_{1}, g_{2}, \dots), h=(h_{1}, h_{2}, \dots)$ are elements of $G$, meaning that $g_{n}, h_{n}\in f(n)\mathbb{Z}$ for all $n\in\mathbb{N}$, and only finitely many of the $g_{n}$ and $h_{n}$ are nonzero. Equipped with this metric, $G$ becomes a discrete metric space. Furthermore, $G$ has bounded geometry since the sequence $\{f(n)\}$ is strictly increasing and every element of $G$ has only finitely many nonzero coordinates. Moreover, $G$ has asymptotic dimension growth $f^{-1}$. 
\end{example}

\subsection{The case of coarse spaces}
Our objective in this subsection is to establish the fourth equivalent definition of asymptotic dimension growth. To achieve this, we consider a coarse space whose coarse structure is induced by a metric. Furthermore, we reformulate the two conditions presented in Definition \ref{def 2.6} using the language of coarse geometry, and demonstrate their equivalence. This new definition is applied in the proof of Theorem \ref{th 4.17}.

We initially recollect the notion of coarse spaces, referring the reader to the monograph \cite{Roe03} for a comprehensive treatment of coarse spaces.
\begin{definition}[\cite{Roe03}]
Let $X$ be a set, and let $\mathcal{E}$ be a collection of subsets of $X\times X$. We say that $\mathcal{E}$ is a coarse structure on $X$ if it satisfies the following axioms:
\begin{itemize}
\item The diagonal $\bigtriangleup=\{(x,x)\mid x\in X\}$ is a member of $\mathcal{E}$;
\item $E\in\mathcal{E}$ and $F\subseteq E$ implies $F\in\mathcal{E}$;
\item $E,F\in\mathcal{E}$ implies $E\cup F\in \mathcal{E}$;
\item $E\in\mathcal{E}$ implies $E^{-1}:=\{(y,x)\mid (x,y)\in E\}\in\mathcal{E}$;
\item $E, F\in\mathcal{E}$ implies $E\circ F:=\{(x,y)\mid \exists z, s.t. (x,z)\in E \text{ and } (z,y)\in F\}\in\mathcal{E}$;
\end{itemize}
The pair $(X,\mathcal{E})$ is called a coarse space, and the elements of $\mathcal{E}$ are called controlled sets or entourages.
\end{definition}

\begin{remark}\label{rem 2.11}
Below, we examine the special case of the coarse structure derived from a metric.
Let $(X,d)$ be a discrete metric space. For $R>0$, we define a controlled tube as follows:
$$E_{R}:=\{(x,y)\in X\times X\mid d(x,y)\leq R\}.$$
Then, the controlled sets are defined as subsets $E$ of $X\times X$ contained within some controlled tube $E_{R}$ for some $R$, thereby forming a coarse structure on $X$.
\end{remark}

The following definition introduces two concepts for a family of subsets in coarse geometry. These are $E$-separatedness and $E$-boundedness.
\begin{definition}[\cite{Roe03}]
Let $X$ be a coarse space, and let $E$ be a controlled set for $X$. A family $\mathcal{U}=\{U_{i}\}_{i\in I}$ of subsets of $X$ is:
\begin{enumerate}
\item $E$-separated if $U_{i}\times U_{j}\cap E=\varnothing$ when $i\neq j$;
\item $E$-bounded if $U_{i}\times U_{i}\subseteq E$ for each $i\in I$.
\end{enumerate}
\end{definition}

The fourth equivalent definition translates the concepts of $R$-disjointness and uniform boundedness from Definition \ref{def 2.6} into the terms of $E_{R}$-separatedness and $F$-boundedness within the framework of coarse geometry.

\begin{definition}\label{def 2.15}
Let $(X, d)$ be a separable metric space. For any $R>0$, we define $\widecheck{\ad}_{X}(R)$ to be the smallest number $m\in\mathbb{N}$ with the following property: for the controlled set $E_{R}$, there exist a controlled set $F$ and a cover $\mathcal{U}=\{U_{i}\}_{i\in I}$ of $X$ such that $\mathcal{U}$ is $F$-bounded and such that $\mathcal{U}$ admits a decomposition
$$
\mathcal{U}=\mathcal{U}_{0}\sqcup\cdots\sqcup\mathcal{U}_{m}
$$
and each $\mathcal{U}_{i}$ is $E_{R}$-separated. 
\end{definition}

The following lemma confirms that this coarse geometric formulation is equivalent to the previous three definitions.

\begin{lemma}\label{lemma 2.16}
Let $(X, d)$ be a separable metric space. Then $\widecheck{\ad}_{X}= \widehat{\ad}_{X}$, hence $\widecheck{\ad}_{X}\approx\widetilde{\ad}_{X}\approx\ad_{X}$.
\end{lemma}
\begin{proof}
Given $R>0$. Suppose that $\mathcal{U}_{0},\dots,\mathcal{U}_{\widehat{\ad}_{X}(R)}$ are uniformly bounded, $R$-disjoint families of subsets of $X$ covering $X$. Since each family $\mathcal{U}_{i}$ is uniformly bounded, the cover $\mathcal{U}=\bigsqcup\limits_{i=0}^{\widehat{\ad}_{X}(R)}\mathcal{U}_{i}$ is also uniformly bounded. Consequently, there exists a constant $C>0$ such that for each $U\in\mathcal{U}$, $\diam(U)\leq C$, which implies that each set $U\times U$ is contained in a controlled set $E_{C}$, thus the cover $\mathcal{U}$ is $E_{C}$-bounded. Furthermore, as each family $\mathcal{U}_{i}$ is $R$-disjoint, then for every $U\neq U'$ in $\mathcal{U}_{i}$, $d(U,U')>R$. This ensures that $(U\times U')\cap E_{R}=\varnothing$, indicating that each $\mathcal{U}_{i}$ is $E_{R}$-separated. Hence, $\widecheck{\ad}_{X}(R)\leq\widehat{\ad}_{X}(R)$, which yields $\widecheck{\ad}_{X}\leq \widehat{\ad}_{X}$.

Conversely, given $R>0$, for the controlled set $E_{R}$, there exist a controlled set $F$ and a cover $\mathcal{U}=\{U_{i}\}_{i\in I}$ of $X$ such that $\mathcal{U}$ is $F$-bounded and admits a decomposition
$$
\mathcal{U}=\mathcal{U}_{0}\sqcup\cdots\sqcup\mathcal{U}_{\widecheck{\ad}_{X}(R)}
$$
such that each $\mathcal{U}_{i}$ is $E_{R}$-separated. Since $F$ is a controlled set, $F\subseteq E_{C}$ for some $C>0$. As $\mathcal{U}$ is $F$-bounded, each set $U_{i}\times U_{i}$ is contained in $ E_{C}$, implying that $\diam(U_{i})\leq C$, so $\mathcal{U}$ is a uniformly bounded cover. Since each $\mathcal{U}_{i}$ is $E_{R}$-separated,  for every $U\neq U'$ in $\mathcal{U}_{i}$, $(U\times U')\cap E_{R}=\varnothing$, thus $d(U, U')>R$, which implies that each $\mathcal{U}_{i}$ is $R$-disjoint.
Therefore, $\widehat{\ad}_{X}(R)\leq \widecheck{\ad}_{X}(R)$, which yields $ \widehat{\ad}_{X}\leq \widecheck{\ad}_{X}$. Hence, $\widehat{\ad}_{X}=\widecheck{\ad}_{X}$, which implies that $\widecheck{\ad}_{X}\approx\widetilde{\ad}_{X}\approx\ad_{X}$ by Proposition \ref{prop 2.7}.
\end{proof}

In light of the preceding equivalences, we can use any one of the functions $\ad_{X}$, $\widetilde{\ad}_{X}$, $\widehat{\ad}_{X}$, or $\widecheck{\ad}_{X}$ to define the asymptotic dimension growth.

\section{Dynamic asymptotic dimension growth for group actions}
In this section, we introduce the dynamic asymptotic dimension function and its associated growth type for actions of discrete groups on compact Hausdorff spaces. This framework generalizes the concept of finite dynamic asymptotic dimension for such actions. Furthermore, by combining examples of groups with asymptotic dimension growth from Section \ref{section 2} with the results established in Lemma \ref{lemma 4.13} and Theorem \ref{th 4.17} of Section \ref{section: dad}, we construct several instructive examples of group actions that exhibit exponential dynamic asymptotic dimension growth. 

We begin by revisiting the concept of length functions on discrete groups, which serves as the foundation for our primary definition.
\begin{definition}
Let $\Gamma$ be a discrete group. A proper length function on $\Gamma$ is a function $\ell: \Gamma\rightarrow [0,\infty)$ satisfying:
\begin{enumerate} 
\item $\ell(\gamma)=0$ if and only if $\gamma=e$;
\item $\ell(\gamma_{1}\gamma_{2})\leq\ell(\gamma_{1})+\ell(\gamma_{2})$ for all $\gamma_{1},\gamma_{2}\in\Gamma$;
\item $\ell(\gamma^{-1})=\ell(\gamma)$ for all $\gamma\in\Gamma$;
\item $\{\gamma\in\Gamma\mid\ell(\gamma)\leq R\}$ is finite for all $R\geq 0$.
\end{enumerate}
\end{definition}

The following defines the dynamic asymptotic dimension growth for actions of discrete groups on compact spaces via proper length functions.
\begin{definition}\label{def 3.2}
Let $\Gamma$ be a discrete group, and let $X$ be a compact Hausdorff space. Let $\Gamma\curvearrowright X$ be an action by homeomorphisms, and let $\ell_{\Gamma}: \Gamma\rightarrow [0,\infty)$ be a proper length function.  For any $R>0$, let $\dad_{\Gamma\curvearrowright X}(R)$ be the smallest number $m\in\mathbb{N}$ with the following property: for the finite subset $E=\{\gamma\in\Gamma\mid\ell_{\Gamma}(\gamma)<R\}$ of $\Gamma$, there exists an open cover $\{U_{0},\dots, U_{m}\}$ of $X$ such that for each $i\in\{0,\dots,m\}$, the set
\begin{align*}
\{\gamma\in\Gamma\mid&\text{ there exist }x\in U_{i}, n\in\mathbb{N}, \text{ and }\gamma_{n},\dots,\gamma_{1}\in E\text{ such that }\gamma=\gamma_{n}\cdots\gamma_{2}\gamma_{1}\\
&\text{ and for all }k\in\{1,\dots,n\}, \gamma_{k}\cdots\gamma_{1}x\in U_{i}\}.
\end{align*}
is finite. We call $\dad_{\Gamma\curvearrowright X}:\mathbb{R}^{+}\rightarrow\mathbb{N}$ the dynamic asymptotic dimension function of $\Gamma\curvearrowright X$, and call the growth type of the function $\dad_{\Gamma\curvearrowright X}$ the dynamic asymptotic dimension growth of $\Gamma\curvearrowright X$.
\end{definition}

\begin{remark}
The following statements hold:
\begin{enumerate}
\item The function $\dad_{\Gamma\curvearrowright X}$ is non-decreasing, and the dynamic asymptotic dimension of the action $\Gamma\curvearrowright X$ is defined as follows:
\begin{equation*}
\dasdim (\Gamma\curvearrowright X)=\lim\limits_{R\rightarrow\infty}\dad_{\Gamma\curvearrowright X}(R).
\end{equation*}
\item $\dasdim(\Gamma\curvearrowright X)\leq d$ for some finite $d$ if and only if $\dad_{\Gamma\curvearrowright X}$ is a constant function and $\dad_{\Gamma\curvearrowright X}\leq d$. 
\item $E=E^{-1}$ and $e\in E$ since $\ell$ is a proper length function.
\end{enumerate}
\end{remark}

\begin{remark}
The dynamic asymptotic dimension has been estimated for actions of discrete groups on compact spaces. For instance, for an action of an infinite virtually cyclic group on a compact Hausdorff space,
the dynamic asymptotic dimension is always one if the action has the marker property \cite{ALSS21}. Moreover, for a free continuous action of a countable discrete group on a zero-dimensional second-countable compact Hausdorff space, the dynamic asymptotic dimension is either infinite or equal to the asymptotic dimension of the group itself \cite{CJMST23}. Additionally, for an infinite countable discrete group $\Gamma$ with asymptotic dimension $d$, the dynamic asymptotic dimension of a free, minimal action of $\Gamma$ on the Cantor set is at most $d$ \cite{GWY17}.
\end{remark}

\begin{remark}
If a discrete group $\Gamma$ has infinite asymptotic dimension, then its corresponding group action has infinite dynamic asymptotic dimension, according to Lemma \ref{lemma 4.13} and [\cite{Bon24}, Proposition 3.3].
\end{remark}

An example of a group action with at most exponential dynamic asymptotic dimension growth is obtained by combining Example \ref{example 2.8} with Lemma \ref{lemma 4.13} and Theorem \ref{th 4.17}.
\begin{example}\label{example 3.6}
Let $\Gamma$ be a finitely generated group with the word metric. Then the canonical action of $\Gamma$ on its Stone-\v{C}ech compactification $\beta\Gamma$ has at most exponential dynamic asymptotic dimension growth. Moreover, the dynamic asymptotic dimension of this action is infinite when $\Gamma$ itself has infinite asymptotic dimension.
\end{example}

A further example of a group action with slow dynamic asymptotic dimension growth emerges from the application of Lemma \ref{lemma 4.13} and Theorem \ref{th 4.17} to Example \ref{example 2.11}.
\begin{example}\label{example 3.7}
Let $f: \mathbb{N}\rightarrow\mathbb{R}^{+}$ be a strictly increasing function such that $f\succeq x^{\alpha}$ for some $\alpha>1$. Let $G=\bigoplus\limits_{n\in\mathbb{N}}f(n)\mathbb{Z}$, and endow $G$ with the direct sum metric from \eqref{eq 2.2}. Then the canonical action of $G$ on its Stone-\v{C}ech compactification $\beta G$ has slow dynamic asymptotic dimension growth. Moreover the dynamic asymptotic dimension of this action is infinite since $G$ has infinite asymptotic dimension.
\end{example}

\section{Dynamic asymptotic dimension growth for groupoids} \label{section: dad}
In this section, we formulate the definition of dynamic asymptotic dimension growth for \'etale groupoids, which generalizes the notion previously introduced for actions of discrete groups on compact spaces. Furthermore, for a discrete metric space $X$ of bounded geometry, we prove that the asymptotic dimension growth of $X$ is equivalent to the dynamic asymptotic dimension growth of its coarse groupoid $G(X)$. Consequently, $G(X)$ is amenable when this growth is subexponential.

Recall from \cite{Ren80} that a groupoid consists of a set $\mathcal{G}$ of morphisms, a subset $\mathcal{G}^{(0)}\subseteq\mathcal{G}$ of objects (called the unit space), together with source and range maps $s,r:\mathcal{G}\rightarrow\mathcal{G}^{(0)}$, a composition map
$$
\mathcal{G}^{(2)}:=\{(g_{1},g_{2})\in\mathcal{G}\times\mathcal{G}\mid s(g_{1})=r(g_{2})\}\ni (g_{1},g_{2})\mapsto g_{1}g_{2}\in\mathcal{G},
$$
and an inverse map $g\mapsto g^{-1}$ satisfying a couple of axioms. Throughout this paper, we assume that every groupoid is locally compact and Hausdorff, meaning that it is equipped with a locally compact and Hausdorff topology such that structure maps are continuous. 

A groupoid $\mathcal{G}$ is called \'etale if the range map $r$ (and consequently the source map $s$) is a local homeomorphism. More precisely, for every $g\in\mathcal{G}$, there exists a neighborhood $U$ of $g$ such that $r(U)$ is open in $\mathcal{G}^{(0)}$ and the restriction $r|_{U}$ is a homeomorphism. For a locally compact Hausdorff groupoid $\mathcal{G}$, it is \'etale if and only if its topology admits a basis consisting of open bisections. We say that $\mathcal{G}$ is ample if its unit space $\mathcal{G}^{(0)}$ is totally disconnected (for example, a Cantor set or a discrete space). For each $x\in\mathcal{G}^{(0)}$, we define $\mathcal{G}^{x}:=r^{-1}(x)$ (the range fiber at $x$) and $\mathcal{G}_{x}:=s^{-1}(x)$ (the source fiber at $x$). Their intersection $\mathcal{G}_{x}^{x}:=\mathcal{G}^{x}\cap \mathcal{G}_{x}=\{g\in\mathcal{G}\mid s(g)=r(g)=x\}$ is called the isotropy group at $x$. The groupoid $\mathcal{G}$ is said to be principal if for every $x\in\mathcal{G}^{(0)}$, the isotropy group $\mathcal{G}_{x}^{x}$ is trivial. 

\subsection{The case of \'etale groupoids}
We now focus on \'etale groupoids, for which it is easier to find sufficient conditions for the amenability of groupoids. In this subsection, we introduce the dynamic asymptotic dimension growth for \'etale groupoids. We first define the dynamic asymptotic dimension growth for \'etale groupoids with compact unit space by using length functions on these groupoids. Further, borrowing the notion of Alexandrov groupoid defined in \cite{KKLRU21}, we extend the definition to the non-compact case. Finally, we verify that this notion is a generalization of the dynamic asymptotic dimension growth for actions of discrete groups on compact spaces.

We now recall the definition of length functions on groupoids from [\cite{MW26}, Definition 4.3], which generalizes and refines earlier versions for locally compact groupoids [\cite{Tu99}, Definition 3.1] and \'etale groupoids [\cite{OY19}, Definition 2.22].

\begin{definition}[\cite{MW26}]
Let $\mathcal{G}$ be a groupoid. A length function on $\mathcal{G}$ is a function $\ell:\mathcal{G}\rightarrow [0,\infty)$ satisfying
\begin{enumerate} 
\item $\ell(x)=0$ if and only if $x\in\mathcal{G}^{(0)}$;
\item $\ell(x)=\ell(x^{-1})$ for all $x\in\mathcal{G}$;
\item $\ell(xy)\leq\ell(x)+\ell(y)$ if $x$ and $y$ are composable in $\mathcal{G}$.
\end{enumerate} 
\end{definition}

The length functions on \'etale groupoids enjoy the following properties.
\begin{definition}[\cite{MW26}]
Let $\mathcal{G}$ be an \'etale groupoid. Let $\ell:\mathcal{G}\rightarrow [0,\infty)$ be a length function on $\mathcal{G}$. For any subset $K\subseteq\mathcal{G}$, we write
$$
\bar{\ell}(K)=\sup\limits_{x\in K}\ell(x).
$$
We say $\ell$ is
\begin{enumerate}
\item proper if for any $K\subseteq\mathcal{G}\backslash\mathcal{G}^{(0)}$, $\bar{\ell}(K)<\infty$ implies that $K$ is precompact;
\item controlled if for any $K\subseteq\mathcal{G}$, $K$ is precompact implies that $\bar{\ell}(K)<\infty$;
\item coarse if it is both proper and controlled;
\item continuous if it is a continuous function with regard to the topology of $\mathcal{G}$.
\end{enumerate}
\end{definition}

\begin{remark}
Recall from \cite{MW26} that two length functions $\ell_{1}, \ell_{2}$ are said to be coarsely equivalent if for every $r>0$, 
$$
\sup\{\ell_{1}(x)\mid x\in\mathcal{G}, \ell_{2}(x)\leq r\}<\infty \text{ and }\sup\{\ell_{2}(x)\mid x\in\mathcal{G}, \ell_{1}(x)\leq r\}<\infty. 
$$ 
When the unit space $\mathcal{G}^{(0)}$ is compact, the item (1) simplifies to: 
$\ell$ is proper if for any $K\subseteq\mathcal{G}$, $\bar{\ell}(K)<\infty$ implies that $K$ is precompact. Note that a continuous length function is controlled, hence the item (4) implies the item (2). 
\end{remark}

In studying the amenability for groupoids, coarse continuous length functions play a fundamental role. Their existence for a certain class of groupoids was established by Ma and Wu in \cite{MW26}:
\begin{theorem}[\cite{MW26}]
Up to coarse equivalence, any $\sigma$-compact locally compact Hausdorff \'etale groupoid has a unique coarse continuous length function.
\end{theorem}

The following transformation groupoid is a typical example of the above theorem. This groupoid plays an important role in the proof of Lemma \ref{lemma 4.13}. 
\begin{example}[\cite{MW26}]\label{example 4.5}
Let $\Gamma$ be a countable discrete group, and let $X$ be a $\sigma$-compact, locally compact Hausdorff space. The action $\Gamma\curvearrowright X$ by homeomorphisms induces an associated transformation groupoid $X\rtimes\Gamma$ as follows:
\begin{itemize}
\item $X\rtimes\Gamma:=\{(\gamma x,\gamma,x)\mid x\in X,\gamma\in\Gamma\}$.
\item The unit space ${(X\rtimes\Gamma)}^{(0)}$ is $\{(x,e,x)\mid x\in X\}$, which is identified with $X$.
\item Source and range maps: $s(\gamma x,\gamma,x)=x$, $r(\gamma x,\gamma,x)=\gamma x$.
\item Composition: $(\gamma\gamma'x,\gamma,\gamma'x)(\gamma'x,\gamma',x)=(\gamma\gamma'x,\gamma\gamma',x)$.
\item Inverse: $(\gamma x,\gamma,x)^{-1}=(x,\gamma^{-1},\gamma x)$.
\end{itemize}
Equipped with the topology induced by $X\times\Gamma$, the groupoid $X\rtimes\Gamma$ becomes a $\sigma$-compact, locally compact, Hausdorff, and \'etale groupoid. Let $\ell_{\Gamma}:\Gamma\rightarrow[0,\infty)$ be a proper length function on $\Gamma$, and let $g: X\rightarrow [0,\infty)$ be a continuous proper function on $X$. We define
\begin{align*}
\ell_{X\rtimes\Gamma}:X\rtimes\Gamma&\rightarrow [0,\infty),\\
(\gamma x,\gamma,x)&\mapsto \ell_{\Gamma}(\gamma)\cdot\left(1+\max\{g(\gamma x),g(x)\}\right).
\end{align*}
Clearly, it is a coarse continuous length function on $X\rtimes\Gamma$. If $X$ is compact, we may simply choose $g=0$, and thus
$$
\ell_{X\rtimes\Gamma}(\gamma x,\gamma,x)=\ell_{\Gamma}(\gamma)\quad \text{for all }(\gamma x,\gamma,x)\in X\rtimes\Gamma.
$$
\end{example}

By using coarse continuous length functions on \'etale groupoids with compact unit space, we define the dynamic asymptotic dimension growth for these groupoids.
\begin{definition}\label{def 4.6}
Let $\mathcal{G}$ be an \'etale groupoid with compact unit space $\mathcal{G}^{(0)}$, and let $\ell:\mathcal{G}\rightarrow [0,\infty)$ be a coarse continuous length function. For any $R>0$, let $\dad_{\mathcal{G}}(R)$ be the smallest number $m\in\mathbb{N}$ with the following property: for the open precompact subset $K={\{z\in\mathcal{G}\mid\ell(z)< R\}}$ of $\mathcal{G}$, there exists an open cover $\{U_{0},\dots, U_{m}\}$ of $\mathcal{G}^{(0)}$ such that for each $i$, the set $\{g\in K\mid s(g), r(g)\in U_{i}\}$ is contained in a precompact subgroupoid of $\mathcal{G}$. We call $\dad_{\mathcal{G}}:\mathbb{R}^{+}\rightarrow\mathbb{N}$ the dynamic asymptotic dimension function of $\mathcal{G}$, and call the growth type of the function $\dad_{\mathcal{G}}$ the dynamic asymptotic dimension growth of $\mathcal{G}$. 
\end{definition}

The following definition from [\cite{KKLRU21}, Definition 7.8] provides a method to handle the case of an \'etale groupoid $\mathcal{G}$ with non-compact unit space by reducing the unit space to the compact case. 

\begin{definition}[\cite{KKLRU21}]
Let $\mathcal{G}$ be an \'etale groupoid with locally compact Hausdorff space of units. Then the Alexandrov groupoid $\mathcal{G}^{+}$ is the set $\mathcal{G}\cup(\mathcal{G}^{(0)})^{+}$ with the topology determined by specifying that the inclusions $\mathcal{G}$, $(\mathcal{G}^{(0)})^{+}\subseteq\mathcal{G}^{+}$ are open, and the groupoid structure extending the groupoid structure of $\mathcal{G}$ and making $\infty$ a unit.
\end{definition}

Observe that $\mathcal{G}^{+}$ is an \'etale groupoid whose unit space is the compact Hausdorff space $(\mathcal{G}^{(0)})^{+}$, obtained by taking Alexandrov one-point compactification of $\mathcal{G}^{(0)}$. 

The following example from [\cite{CDGHV24}, Example 3.5] illustrates that for the transformation groupoid $\mathcal{G}=X\rtimes\Gamma$ with $X$ non-compact, its Alexandrov groupoid  $\mathcal{G}^{+}$ and the transformation groupoid $X_{+}\rtimes\Gamma$ have the same compact unit space $X_{+}$, but the two groupoids are not identical.

\begin{example}[\cite{CDGHV24}]\label{example 4.8}
Let $\Gamma$ be a discrete group with identity $e$ which acts continuously on a locally compact, non-compact, Hausdorff space $X$ and let $\mathcal{G}=X\rtimes \Gamma$. Let $\mathcal{G}^{+}$ be the Alexandrov groupoid of $\mathcal{G}$, and let $X_{+}$ be the one-point compactification of $X$. It is easy to see that the action of $\Gamma$ on $X_{+}$ given by 
$$
\gamma\cdot x=\begin{cases}\gamma\cdot x& \text{ if }x\in X\\
\infty&\text{ if }x=\infty
\end{cases}
$$
is continuous. Then we have
$$
\mathcal{G}^{+}=\{(\gamma,x)\mid \gamma\in\Gamma, x\in X\}\cup\{(e,\infty)\}\text{ and }X_{+}\rtimes\Gamma=\{(\gamma,x)\mid \gamma\in\Gamma, x\in X_{+}\}.
$$
It is straightforward to verify that the Alexandrov groupoid $\mathcal{G}^{+}$ is not identical to the transformation groupoid $X_{+}\rtimes\Gamma$; it is merely its wide and open subgroupoid.
\end{example}

The following theorem, which relies on [\cite{CDGHV24}, Theorem 3.1], demonstrates that the original groupoid $\mathcal{G}$ shares the same dynamic asymptotic dimension with its Alexandrov groupoid $\mathcal{G}^{+}$.
\begin{theorem}[\cite{CDGHV24}]\label{th 4.9}
Let $\mathcal{G}$ be a locally compact, Hausdorff and \'etale groupoid with non-compact unit space $\mathcal{G}^{(0)}$. Then $\dad(\mathcal{G})=\dad(\mathcal{G}^{+})$.
\end{theorem}

To define the dynamic asymptotic dimension growth for a groupoid $\mathcal{G}$ with non-compact unit space, we apply Theorem \ref{th 4.9} to Definition \ref{def 4.6} by setting $\dad_{\mathcal{G}}(R):=\dad_{\mathcal{G}^{+}}(R)$ for all $R>0$, thereby extending the definition via the Alexandrov groupoid $\mathcal{G}^{+}$.

\begin{definition}\label{def 4.10}
Let $\mathcal{G}$ be an \'etale groupoid with non-compact unit space $\mathcal{G}^{(0)}$, and let $\mathcal{G}^{+}$ be the Alexandrov groupoid of $\mathcal{G}$. We define the dynamic asymptotic dimension function of $\mathcal{G}$ as $\dad_\mathcal{G}:=\dad_{\mathcal{G}^{+}}$, and call the growth type of the function $\dad_{\mathcal{G}}$ the dynamic asymptotic dimension growth of $\mathcal{G}$.
\end{definition}

\begin{remark}
The following statements hold:
\begin{enumerate}
\item The function $\dad_{\mathcal{G}}$ is non-decreasing, and the dynamic asymptotic dimension of $\mathcal{G}$ is defined as follows:
\begin{equation*}
\dasdim(\mathcal{G})=\lim\limits_{R\rightarrow\infty}\dad_{\mathcal{G}}(R).
\end{equation*}
\item $\dasdim(\mathcal{G})\leq d$ for some finite $d$ if and only if $\dad_{\mathcal{G}}$ is a constant function and $\dad_{\mathcal{G}}\leq d$. 
\item For each $i$, the subgroupoid $\mathcal{G}_{i}=\langle K\cap\mathcal{G}|_{U_{i}}\rangle$ generated by the set $$\{g\in K\mid s(g), r(g)\in U_{i}\}$$ is open precompact by the next lemma.
\end{enumerate} 
\end{remark}

The following lemma is used in the proof of Lemma \ref{lemma 5.6}. 
\begin{lemma}[\cite{GWY17}]\label{lemma 4.11}
Let $\mathcal{G}$ be an \'etale groupoid.
\begin{enumerate}
\item If $K, H$ are open (resp. compact, precompact) subsets of $\mathcal{G}$, then
$$K\cdot H:=\{kh\in\mathcal{G}\mid k\in K, h\in H, s(k)=r(h)\}$$
is open (resp. compact, precompact).
\item If $K\subseteq\mathcal{G}$ is open, then the subgroupoid of $\mathcal{G}$ generated by $K$ is also open, and is itself an \'etale groupoid.
\end{enumerate} 
\end{lemma}

Using the length function on the transformation groupoid constructed in Example \ref{example 4.5}, we prove that the group action $\Gamma\curvearrowright X$ shares the same dynamic asymptotic dimension function with its corresponding transformation groupoid $X\rtimes\Gamma$. This means that Definition \ref{def 4.6} generalizes the dynamic asymptotic dimension growth for group actions in the sense of Definition \ref{def 3.2}.

\begin{lemma}\label{lemma 4.13}
Let $X$ be a compact Hausdorff space, and let $f:\mathbb{R}^{+}\rightarrow\mathbb{N}$ be a nondecreasing function. An action $\Gamma\curvearrowright X$ has dynamic asymptotic dimension function $f$ if and only if the corresponding transformation groupoid $\mathcal{G}=X\rtimes\Gamma$ has dynamic asymptotic dimension function $f$.
\end{lemma}

\begin{proof}
We may assume that an action $\Gamma\curvearrowright X$ has dynamic asymptotic dimension function $\dad_{\Gamma\curvearrowright X}$, and its associated transformation groupoid has dynamic asymptotic dimension function $\dad_{\mathcal{G}}$. It thus suffices to show that $\dad_{\mathcal{G}}=\dad_{\Gamma\curvearrowright X}$. First, we prove that $\dad_{\mathcal{G}}\leq \dad_{\Gamma\curvearrowright X}$. For any $R>0$, for the open precompact subset $K={\{(\gamma x,\gamma,x)\in \mathcal{G}\mid\ell(\gamma x,\gamma,x)< R\}}$ of $\mathcal{G}$. Then $\ell_{\Gamma}(\gamma)=\ell(\gamma x,\gamma,x)<R$, thus $K$ has the form
$$
\{(\gamma x,\gamma,x)\in \mathcal{G}\mid x\in X, \gamma\in E\}
$$ 
for a finite subset $E=\{\gamma\in\Gamma\mid\ell_{\Gamma}(\gamma)<R\}$ of $\Gamma$. Since the action has dynamic asymptotic dimension function $\dad_{\Gamma\curvearrowright X}$, there is an open cover $U_{0},\dots,U_{\dad_{\Gamma\curvearrowright X}(R)}$ of $\mathcal{G}^{(0)}=X$ such that the set
\begin{align*}
F_{i}:=\{\gamma\in\Gamma\mid&\text{ there exist }x\in U_{i}, n\in\mathbb{N}, \text{ and }\gamma_{n},\dots,\gamma_{1}\in E\\&\text{ such that }\gamma=\gamma_{n}\cdots\gamma_{2}\gamma_{1}
\text{ and for all }k\in\{1,\dots,n\},\\&\gamma_{k}\cdots\gamma_{1}x\in U_{i}\}
\end{align*}
is finite for each $i\in\{0,\cdots,\dad_{\Gamma\curvearrowright X}(R)\}$.
Therefore, the subgroupoid $\mathcal{G}_{i}$ generated by the set 
$$
\{(\gamma x,\gamma,x)\in K\mid x,\gamma x\in U_{i}\}
$$
is contained in $\{(\gamma x,\gamma,x)\in\mathcal{G}\mid x\in X, \gamma\in F_{i}\}.$
Since $F_{i}$ is finite and $X$ is compact, the set $\{(\gamma x,\gamma,x)\in\mathcal{G}\mid x\in X, \gamma\in F_{i}\}$ is precompact, and thus $\mathcal{G}_{i}$ is precompact. Hence, $\dad_{\mathcal{G}}(R)\leq \dad_{\Gamma\curvearrowright X}(R)$, which implies that $\dad_{\mathcal{G}}\leq \dad_{\Gamma\curvearrowright X}$.

On the other hand, we show that $\dad_{\Gamma\curvearrowright X}\leq\dad_{\mathcal{G}}$. For any $R>0$, for the finite subset $E=\{\gamma\in\Gamma\mid\ell_{\Gamma}(\gamma)<R\}$ of $\Gamma$, there is an open precompact subset $K=\{(\gamma x,\gamma,x)\in\mathcal{G}\mid x\in X, \gamma\in E\}$ of $\mathcal{G}$. Since $\mathcal{G}$ has the dynamic asymptotic dimension function $\dad_{\mathcal{G}}$, there is an open cover $\{U_{0},\dots, U_{\dad_{\mathcal{G}}(R)}\}$ of $X$ such that for each $i\in\{0,\dots,\dad_{\mathcal{G}}(R)\}$, the set
$$\{(\gamma x,\gamma,x)\in K\mid x,\gamma x\in U_{i}\}$$
is contained in a precompact subgroupoid of $\mathcal{G}$. Then the subgroupoid $\mathcal{G}_{i}$ generated by the set
$$
\{(\gamma x,\gamma,x)\in\mathcal{G}\mid\gamma\in E, x, \gamma x\in U_{i}\}, 
$$
is precompact. For each $i$, let 
\begin{align*}
F_{i}:=\{\gamma\in\Gamma\mid&\text{ there exist }x\in U_{i}, n\in\mathbb{N}, \text{ and }\gamma_{n},\dots,\gamma_{1}\in E\text{ such that }\gamma=\gamma_{n}\cdots\gamma_{2}\gamma_{1}\\
&\text{ and for all }k\in\{1,\dots,n\}, \gamma_{k}\cdots\gamma_{1}x\in U_{i}\}.
\end{align*}
Therefore, the subset $$
\{(\gamma x,\gamma,x)\in\mathcal{G}\mid\gamma\in F_{i}, x\in U_{i}\}
$$
is contained in $\mathcal{G}_{i}$, and is thus precompact, which implies that $F_{i}$ is finite. Clearly, $\dad_{\Gamma\curvearrowright X}(R)\leq \dad_{\mathcal{G}}(R)$, which forces $\dad_{\Gamma\curvearrowright X}\leq\dad_{\mathcal{G}}$. Hence, $\dad_{\Gamma\curvearrowright X}=\dad_{\mathcal{G}}$, namely, $\dad_{\Gamma\curvearrowright X}=f$ if and only if $\dad_{\mathcal{G}}=f$.
\end{proof}

\begin{remark}
For the transformation groupoid $\mathcal{G}=X\rtimes\Gamma$ with $X$ non-compact, we define its dynamic asymptotic dimension growth using the Alexandrov groupoid $\mathcal{G}^{+}$ in Definition \ref{def 4.10}. However, $\mathcal{G}^{+}$ itself is not a transformation groupoid, as shown in Example \ref{example 4.8}. Consequently, Definition \ref{def 4.10} does not yield a definition of dynamic asymptotic dimension growth for the action $\Gamma\curvearrowright X$ on the non-compact Hausdorff space $X$, and thus this approach cannot be used to directly extend Lemma \ref{lemma 4.13} to the non-compact case. Moreover, in the compact case, the proof of Lemma \ref{lemma 4.13} relies on the equality of the length functions of the group and the transformation groupoid. However, when $X$ is non-compact, this equality fails because the length function of $X\rtimes\Gamma$ becomes a bivariate function depending on both the length function of $\Gamma$ and a continuous proper function $g$ on $X$ (see Example \ref{example 4.5}). Hence, the current proof of this lemma fails in the non-compact case.
\end{remark}

\subsection{Example: coarse spaces with asymptotic dimension growth}
In this subsection, we study a geometric class of groupoids. For a discrete metric space $X$ of bounded geometry, Skandalis, Tu, and Yu introduced an \'etale groupoid $G(X)$ in \cite{STY02}, called the coarse groupoid since it encodes many coarse geometric properties of the underlying space $X$. Moreover, we show that the asymptotic dimension growth for $X$ is equivalent to the dynamic asymptotic dimension growth for its associated coarse groupoid $G(X)$. This example not only serves as our primary motivation but also inspires the terminology `dynamic asymptotic dimension growth'.

\begin{definition}[\cite{STY02}]
Let $(X,d)$ be a discrete metric space of bounded geometry. The coarse groupoid $G(X)$ associated to $(X,d)$ is defined as follows:
\begin{enumerate}
\item for any $R>0$, we define $E_{R}:=\{(x,y)\in X\times X\mid d(x,y)\leq R\}$;
\item as a topological space, we have $G(X):=\bigcup\limits_{R\geq 0}\overline{E_{R}}$ in $\beta(X\times X)$, the Stone-\v{C}ech compactification of $X\times X$;
\item equipping $G(X)$ with the weak topology coming from this union yields a locally compact topology for which each $\overline{E_{R}}$ is a compact open subset of $G(X)$;
\item $G(X)^{(0)}=\overline{E_{0}}\cong\beta X$;
\item $r(x,y)=x,s(x,y)=y$, for all $(x,y)\in G(X)$;
\item $E_{R}\circ E_{S}\subseteq E_{R+S}$, for all $R, S\geq 0$.
\end{enumerate}
\end{definition}

\begin{remark}
The coarse groupoid $G(X)$, as constructed in \cite{STY02}, is a $\sigma$-compact, locally compact, Hausdorff, principal, ample, and \'etale groupoid. Here a discrete metric space $X$ of bounded geometry is also called a uniformly locally finite metric space. Recall from \cite{MW26} that there is a canonical length function $\ell$ on $G(X)$ defined by extending the metric $d: X\times X\rightarrow [0,\infty]$ to $\beta(X\times X)$ and observing that it takes finite values on $G(X)$. The canonical length function $\ell$ on $G(X)$ has the following properties:
 \begin{enumerate}
\item $\ell$ is continuous since it comes from the metric $d$.
\item $\ell$ is coarse. Since $X\times X$ is dense in $\beta(X\times X)$, for any $R>0$, the open set $\ell^{-1}([0,R))$ is contained in $\overline{\ell^{-1}([0,R))\cap X\times X}$, which is contained in the compact set $\overline{E_{R}}$. Hence, $\ell^{-1}([0,R))$ is precompact, which yields that $\ell$ is proper. Since $\ell$ is continuous, it is controlled. 
\item $\ell^{-1}([0,R])=\overline{E_{R}}$ and $\ell^{-1}(\{0\})=\overline{E}_{0}\cong\beta X$.
\end{enumerate}
\end{remark}

The following theorem bridges the dynamic asymptotic dimension growth for the coarse groupoid $G(X)$ and the asymptotic dimension growth for the underlying space $X$, which generalizes the finite-dimensional case of [\cite{GWY17}, Theorem 6.4].
\begin{theorem}\label{th 4.17}
Let $X$ be a discrete metric space of bounded geometry, and let $G(X)$ be the associated coarse groupoid. Let $f: \mathbb{R}^{+}\rightarrow\mathbb{N}$ be a non-decreasing function. Then $\ad_{X}\approx f$ if and only if $\dad_{G(X)}\approx f$. 
\end{theorem}

\begin{remark}
For a discrete metric space $X$ of bounded geometry with at most exponential asymptotic dimension growth and infinite asymptotic dimension, the corresponding coarse groupoid $G(X)$ exhibits at most exponential dynamic asymptotic dimension growth and infinite dynamic asymptotic dimension, respectively, as guaranteed by Theorem \ref{th 4.17} and [\cite{Bon24}, Proposition 3.3].
\end{remark}

\begin{proof}
Let us assume that $G(X)$ has dynamic asymptotic dimension function $\dad_{G(X)}$, and $X$ has asymptotic dimension function $\widecheck{\ad}_{X}$ as defined in Definition \ref{def 2.15}. Below, we show that $\widecheck{\ad}_{X}\approx\dad_{G(X)}$. The first step is to establish the inequality $\widecheck{\ad}_{X}\preceq\dad_{G(X)}$. For any $R>0$, for the open precompact subset $K=\{z\in G(X)\mid\ell(z)< R\}$ of $G(X)$, there exists a controlled set $E_{R-\frac{1}{n}}$ for some positive integer $n$ such that $\overline{E_{R-\frac{1}{n}}}\subseteq K$. Observe that $\overline{E_{R-\frac{1}{n}}}$ is a compact open subset of $G(X)$ containing the unit space $G(X)^{(0)}=\beta X$. Let $U_{0},\dots,U_{\dad_{G(X)}(R)}$ be an open cover of $G(X)^{(0)}=\beta X$ such that for each $i\in\{0,\dots,\dad_{G(X)}(R)\}$, the set
$$
\{g\in K\mid s(g),r(g)\in U_{i}\}
$$
is contained in a precompact subgroupoid of $G(X)$. Therefore, the subgroupoid $\mathcal{G}_{i}$ generated by the set
$$
\{g\in \overline{E_{R-\frac{1}{n}}}\mid s(g),r(g)\in U_{i}\}
$$
is precompact. By the definition of the topology on $G(X)$, there exists a controlled set $F\subseteq X\times X$ such that each $\mathcal{G}_{i}$ is contained in $\overline{F}\subseteq G(X)$. Fix $i$, and let $\sim$ denote the equivalence relation on $U_{i}$ induced by $\mathcal{G}_{i}$, so $x\sim y$ if there exists $g\in \mathcal{G}_{i}$ such that $s(g)=x, r(g)=y$. We denote by $\{U^{j}_{i}\mid j\in J_{i}\}$ the equivalence classes for this relation. Define 
$$
\mathcal{U}_{i}:=\{U^{j}_{i}\cap X\mid j\in J_{i}\}.
$$
Let $\mathcal{U}=\mathcal{U}_{0}\sqcup\cdots\sqcup\mathcal{U}_{\dad_{G(X)}(R)}$, and it is a cover of $X$ since $\{U_{0},\dots, U_{\dad_{G(X)}(R)}\}$ is a cover of $\beta X$.

We next claim that the cover $\mathcal{U}$ is $F$-bounded and each $\mathcal{U}_{i}$ is $E_{R-\frac{1}{n}}$-separated. Indeed, the equivalence relation on each $U^{j}_{i}\cap X$ induced by $\mathcal{G}_{i}$ is entirely contained in $F$ as $\mathcal{G}_{i}\subseteq\overline{F}$, so $\mathcal{U}$ is $F$-bounded. For each $\mathcal{U}_{i}$, if $(U^{j}_{i}\times U^{k}_{i})\cap E_{R-\frac{1}{n}}\neq\varnothing$ for $j\neq k$, then there exists $(x,y)\in (U^{j}_{i}\times U^{k}_{i})\cap E_{R-\frac{1}{n}}$. Observe that $s(x,y)=y, r(x,y)=x$, so $(x,y)\in \mathcal{G}_{i}$ and $x\sim y$. However, when $j\neq k$, $U^{j}_{i}\cap U^{k}_{i}=\varnothing$, contradiction. Thus, for $j\neq k$, $(U^{j}_{i}\times U^{k}_{i})\cap E_{R}=\varnothing$, i.e. $\mathcal{U}_{i}$ is $E_{R-\frac{1}{n}}$-separated. Hence, $\widecheck{\ad}_{X}(R-\frac{1}{n})\leq \dad_{G(X)}(R)$, which implies that $\widecheck{\ad}_{X}\preceq\dad_{G(X)}$. 

The next step is to establish the converse inequality $\dad_{G(X)}\preceq \widecheck{\ad}_{X}$. For any $R>0$, for the open precompact subset $K=\{z\in G(X)\mid\ell(z)< R\}$ of $G(X)$. The definition of the topology on $G(X)$ implies that $K\subseteq\overline{E_{R}}$ for the controlled set $E_{R}$. Since $X$ has asymptotic dimension function $\widecheck{\ad}_{X}$, there exist a controlled set $F$ and a cover $\mathcal{U}$ of $X$ such that $\mathcal{U}$ is $F$-bounded and such that
$$
\mathcal{U}=\mathcal{U}_{0}\sqcup\cdots\sqcup\mathcal{U}_{\widecheck{\ad}_{X}(R)}
$$
and each $\mathcal{U}_{i}$ is $E_{R}$-separated.
For each $i$, set
\begin{equation*}
U_{i}:=\overline{\bigsqcup\limits_{U\in\mathcal{U}_{i}}U},
\end{equation*}
which is a compact open subset of $\beta X$. Since $\mathcal{U}$ is a cover of $X$, $\{U_{0},\dots,U_{\widecheck{\ad}_{X}(R)}\}$ is a cover of $\beta X$. Furthermore, for each $i$,
$$
\mathcal{G}_{i}:=\bigsqcup\limits_{U\in\mathcal{U}_{i}}U\times U
$$
is a subgroupoid of the pair groupoid that is contained in $F$. By the continuity of the groupoid operations, the subgroupoid of $G(X)$ generated by 
$$
\{g\in E_{R}\mid r(g), s(g)\in U_{i}\}
$$
is contained in $\overline{\mathcal{G}_{i}}$, a compact subgroupoid of $G(X)$ contained in the compact subset $\overline{F}$. Thus the subset $\{g\in K\mid s(g), r(g)\in U_{i}\}$ is contained in a compact open subgroupoid of $G(X)$. Therefore, $\dad_{G(X)}(R)\leq \widecheck{\ad}_{X}(R)$ for every $R>0$, which implies the desired inequality $\dad_{G(X)}\preceq \widecheck{\ad}_{X}$. Combined with the previous result, we conclude that $\dad_{G(X)}\approx \widecheck{\ad}_{X}$.

Finally, by Lemma \ref{lemma 2.16}, we see that $\ad_{X}\approx\widecheck{\ad}_{X}$, which implies that $\ad_{X}\approx\dad_{G(X)}$. Therefore, $\ad_{X}\approx f$ if and only if $\dad_{G(X)}\approx f$. 
\end{proof}

\begin{remark}
Based on the second step of the above proof, we see that if $\ad_{X}\preceq f$, then $\dad_{G(X)}\preceq f$ in the following slightly stronger form: For any $R>0$, for the open precompact subset  $K=\{z\in G(X)\mid\ell(z)<R\}$  of $G(X)$, there exists an open cover $\{U_{0},\dots, U_{\widecheck{\ad}_{X}(R)}\}$ of $G(X)^{(0)}$ such that the subset $\{g\in K\mid s(g), r(g)\in U_{i}\}$ is contained in a compact open subgroupoid of $G(X)$.  
\end{remark}

A result of Ozawa and Oppenheim establishes that subexponential asymptotic dimension growth implies property A.
\begin{theorem}[\cite{Oza12, Opp14}]\label{th 4.20}
Let $X$ be a discrete metric space of bounded geometry. If $X$ has subexponential asymptotic dimension growth, then $X$ has property A.
\end{theorem}

The equivalence between property A for $X$ and the amenability of the corresponding coarse groupoid $G(X)$ is established by [\cite{STY02}, Theorem 5.3]. Combining this with Theorems \ref{th 4.17} and \ref{th 4.20}, we obtain the following corollary. This result serves as a reformulation of Theorem \ref{th 4.20}.
\begin{corollary}\label{cor 4.21}
Let $X$ be a discrete metric space of bounded geometry. If the associated coarse groupoid $G(X)$ has subexponential dynamic asymptotic dimension growth, then $G(X)$ is amenable.
\end{corollary}

The following example is a direct consequence of Example \ref{example 2.9}, Theorem \ref{th 4.17}, and Corollary \ref{cor 4.21}.
\begin{example}
Let $X$ be a geodesic uniformly locally finite coarse median space with finite rank and at most exponential volume growth. Then its corresponding coarse groupoid $G(X)$ has subexponential dynamic asymptotic dimension growth, and it is amenable.
\end{example}

\section{Generalized partitions of unity}\label{section 5}
The goal of this section is to generalize partitions of unity from \cite{GWY17}, termed generalized partitions of unity. This tool not only serves as a partition of unity for groupoids with slow dynamic asymptotic dimension growth, but also provides almost invariant properties and useful estimates. Furthermore, it is essential to demonstrate the amenability of groupoids in the following section.

The key techniques needed in the proof of the amenability of groupoids are summarized in the following proposition, which generalizes [\cite{GWY17}, Proposition 7.1]. 

\begin{proposition}\label{prop 5.1}
Let $\mathcal{G}$ be an \'etale groupoid with compact unit space, and with dynamic asymptotic dimension function $f\preceq x^{\alpha} (0<\alpha<1)$. Then for any $R>0$, for the open precompact subset $K={\{z\in\mathcal{G}\mid\ell(z)< R\}}$ of $\mathcal{G}$ and any $\varepsilon>0$, there exist a constant $c$ depending on parameters $\alpha, R, \varepsilon$, a positive integer $p$ depending only on $\alpha$, and an open cover $\{U_{0},\dots, U_{f([(R+1)^{c}]R+1)}\}$ of $\mathcal{G}^{(0)}$ with the following properties:
\begin{enumerate}
\item For each $i$, the set $$\{g\in K\mid s(g), r(g)\in U_{i}\}$$ is contained in an ( open and ) precompact subgroupoid of $\mathcal{G}$.
\item For all $x\in\mathcal{G}^{(0)}$, the `partial orbit' $s(K\cap r^{-1}(x))$ is completely contained in some $U_{i}$.
\item There exists a collection of continuous functions $\{\phi_{i}:\mathcal{G}^{(0)}\rightarrow [0,1]\}_{i=0}^{f([(R+1)^{c}]R+1)}$ such that $\supp(\phi_{i})\subseteq U_{i}$ for each $i$, and $\sum\limits_{i=0}^{f([(R+1)^{c}]R+1)}\left(\phi_{i}(x)\right)^{p}=1$ for all $x\in\mathcal{G}^{(0)}$, and such that for any $g\in K$,
$$
p\cdot\sum\limits_{i=0}^{f([(R+1)^{c}]R+1)}\left|\phi_{i}(r(g))-\phi_{i}(s(g))\right|<\varepsilon.
$$
\end{enumerate}
\end{proposition}

To understand the above proposition more intuitively, we look at the case of group actions on compact spaces.
\begin{corollary}
Let $\Gamma\curvearrowright X$ be an action with $X$ compact, and with dynamic asymptotic dimension function $f\preceq x^{\alpha} (0<\alpha<1)$. Then for any $R>0$, for the finite subset $E=\{\gamma\in\Gamma\mid \ell(\gamma)<R\}$ of $\Gamma$ and any $\varepsilon>0$,  there exist a constant $c$ depending on parameters $\alpha, R, \varepsilon$, a positive integer $p$ depending only on $\alpha$, and an open cover $\{U_{0},\dots, U_{f([(R+1)^{c}]R+1)}\}$ of $X$ with the following properties: 

\begin{enumerate}
\item For each $i$, the set 
\begin{align*}
\{\gamma\in\Gamma\mid&\text{ there exist }x\in U_{i}, n\in\mathbb{N}, \text{ and }\gamma_{n},\dots,\gamma_{1}\in E\\&\text{ such that }\gamma=\gamma_{n}\cdots\gamma_{2}\gamma_{1}
\text{ and for all }k\in\{1,\dots,n\}, \gamma_{k}\cdots\gamma_{1}x\in U_{i}\}
\end{align*}
is finite.
\item For all $x\in X$, the collection $E\cdot x:=\{\gamma^{-1}x\mid\gamma\in E\}$ is completely contained in some $U_{i}$.
\item There exists a collection of continuous functions $\{\phi_{i}:\mathcal{G}^{(0)}\rightarrow [0,1]\}_{i=0}^{f([(R+1)^{c}]R+1)}$ such that $\supp(\phi_{i})\subseteq U_{i}$ for each $i$, and $\sum\limits_{i=0}^{f([(R+1)^{c}]R+1)}\left(\phi_{i}(x)\right)^{p}=1$ for all $x\in X$, and such that for any $\gamma\in E$,
$$
p\cdot\sum\limits_{i=0}^{f([(R+1)^{c}]R+1)}\left|\phi_{i}(\gamma x)-\phi_{i}(x)\right|<\varepsilon.
$$
\end{enumerate}
\end{corollary}

To derive the main proposition of this section, we first prove the following lemma, which serves as an essential tool for subsequent lemmas and theorems.

\begin{lemma}\label{lemma 5.3}
Let $\mathcal{G}$ be an \'etale groupoid, and let $K$ be a symmetric subset of $\mathcal{G}$ containing $\mathcal{G}^{(0)}$. Let $\mathcal{U}=\{U_{j}\}_{j\in J}$ be an open cover of $\mathcal{G}^{(0)}$. Then for any $g\in K$, $s(g)\in U_{j}$ if and only if $r(g)\in U_{j}$ for some $j\in J$.
\end{lemma}

\begin{proof}
For any $x\in U_{j}\subseteq\mathcal{G}^{(0)}\subseteq K$, $s(x)=r(x)=x$, thus $U_{j}\subseteq K\cap r^{-1}(U_{j})$, which implies that $U_{j}=s(U_{j})\subseteq s(K\cap r^{-1}(U_{j}))$. If $g\in K$ and $s(g)\in U_{j}$ for some $j$, then $r(g)$ must be in $U_{j}$. Otherwise, $g\notin K\cap r^{-1}(U_{j})$, and so $s(g)\notin s(K\cap r^{-1}(U_{j}))\supseteq U_{j}$, hence $s(g)\notin U_{j}$, contradiction. Conversely, assume that $g\in K$ and $r(g)\in U_{j}$ for some $j$. Due to symmetry, $g^{-1}\in K$, $s(g^{-1})=r(g)\in U_{j}$. From the previous steps, we see that $s(g)=r(g^{-1})\in U_{j}$.
\end{proof}

Applying the above conclusion, we have the following lemma, which is an adaptation of [\cite{GWY17}, Lemma 7.3].

\begin{lemma}\label{lemma 5.4}
Let $\mathcal{G}$ be an \'etale groupoid with compact unit space, and with dynamic asymptotic dimension function $f$. Then for any $R>0$, for the open precompact subset $K={\{z\in\mathcal{G}\mid\ell(z)< R\}}$ of $\mathcal{G}$, there is a cover $\{U_{0},\dots, U_{f(R)}\}$ of $\mathcal{G}^{(0)}$ by open precompact subsets of $\mathcal{G}^{(0)}$ such that for each $i$, the set 
$$
\{g\in K\mid s(g), r(g)\in U_{i}\}
$$
generates an open precompact subgroupoid of $\mathcal{G}$, and moreover so that for each $x\in\mathcal{G}^{(0)}$, the set $s(K\cap r^{-1}(x))$ is completely contained in some $U_{i}$.
\end{lemma}

\begin{proof}
Since $\mathcal{G}$ has dynamic asymptotic dimension function $f$ and $\mathcal{G}^{(0)}$ is compact, there is an open cover $\{U_{0},\dots, U_{f(R)}\}$ of $\mathcal{G}^{(0)}$ by precompact subsets of $\mathcal{G}^{(0)}$ such that for each $i$, the set
$$
\{g\in K\mid s(g), r(g)\in U_{i}\}
$$
is contained in a precompact subgroupoid of $\mathcal{G}$, and hence generates an open precompact subgroupoid of $\mathcal{G}$. Observe that $K$ is symmetric as $\ell(x)=\ell(x^{-1})$ for all $x\in\mathcal{G}$. Moreover, for all $x\in\mathcal{G}^{(0)}\subseteq K$, we know that $r(x)=x$, so $x\in K\cap r^{-1}(x)$, thus $K\cap r^{-1}(x)\neq\varnothing$. For any $g\in K\cap r^{-1}(x)$, we have $g\in K$ and $r(g)=x\in U_{i}$ for some $i$. By Lemma \ref{lemma 5.3}, $s(g)\in U_{i}$. Hence, $s(K\cap r^{-1}(x))\subseteq U_{i}$ for some $U_{i}$.
\end{proof}

If the unit space of an \'etale groupoid is compact and Hausdorff, then it is normal. This yields the following key result:
\begin{lemma}\label{lemma 5.5}
Let $\mathcal{G}$ be an \'etale groupoid with compact unit space, and let $V_{0},\dots, V_{n}$ be an open cover of $\mathcal{G}^{(0)}$. Then there exist open subsets $\mathcal{U}_{0},\dots,\mathcal{U}_{n}$ of $\mathcal{G}^{(0)}$ such that $\bigcup\limits_{i=0}^{n}\mathcal{U}_{i}$ covers $\mathcal{G}^{(0)}$ and $\overline{\mathcal{U}_{i}}\subseteq V_{i}$ for each $i$.
\end{lemma}

\begin{proof}
We divide our proof into four steps. First, we show that for any open subset $V$ of $\mathcal{G}^{(0)}$ and any $x\in V$, there exists an open neighborhood $U_{x}$ of $x$ such that $\overline{U_{x}}\subseteq V$. Note that $\mathcal{G}^{(0)}$ is Hausdorff, thus the singleton set $\{x\}$ is closed. Additionally, since $\mathcal{G}^{(0)}$ is compact, it is also normal. Observe that the boundary $\partial V$ is closed, and $\{x\}\cap\partial V\subseteq V\cap\partial V=\varnothing$ as $V$ is open. Hence by normality of $\mathcal{G}^{(0)}$, there exist an open neighborhood $U_{x}$ of $x$ and an open set $V'$ containing $\partial V$ such that 
$$U_{x}\cap V'=\varnothing.$$    
We may assume that $U_{x}\subseteq V$, otherwise we can replace $U_{x}$ with $U_{x}\cap V$. It suffices to show that $\overline{U_{x}}\subseteq V$. If this does not happen, we may find a point $y\in\partial U_{x}\subseteq\overline{V}$ such that $y\notin V$, which implies $y\in\partial V\subseteq V'$. It follows that $V'$ is an open neighborhood of $y$, and that $U_{x}\cap V'\neq\varnothing$ as $y\in\partial U_{x}$, contradiction. Hence, we conclude that $\overline{U_{x}}\subseteq V$, as required. 

In the second step, we build subsets $V_{0}'$ and $\mathcal{U}_{0}$ of $\mathcal{G}^{(0)}$ such that $$V_{0}'\subseteq\mathcal{U}_{0}\subseteq\overline{\mathcal{U}_{0}}\subseteq V_{0}.$$ If $y\in\partial V_{0}$, namely $y\notin V_{0}$, then $y$ must be in $V_{i}$ for some $i\in\{1,\dots,n\}$. Thus there exists an open neighborhood $U_{y}'$ of $y$ such that $U_{y}'\subseteq V_{i}$ as $V_{i}$ is open. Observe that for any $z\in\partial V_{i}$, $z$ and $y$ are distinct points. Since $\mathcal{G}^{(0)}$ is Hausdorff, there exist open neighborhoods $\widetilde{U_{z}}$ of $z$ and $\widetilde{U_{y}}$ of $y$ such that $\widetilde{U_{z}}\cap\widetilde{U_{y}}=\varnothing$. Set $U_{y}=\widetilde{U_{y}}\cap U_{y}'$, then $U_{y}\cap\widetilde{U_{z}}=\varnothing$ and $U_{y}\subseteq V_{i}$. Let $V_{0}'=\overline{V_{0}}\backslash\bigcup\limits_{y\in\partial V_{0}}U_{y}$, it is closed and hence compact. Note that $V_{0}'\subseteq V_{0}$, then by the first step, for any $x\in V_{0}'$, there exists an open neighborhood $U_{x}$ of $x$ such that $\overline{U_{x}}\subseteq V_{0}$ and $\bigcup\limits_{x\in V_{0}'}U_{x}$ is a cover of $V_{0}'$. In view of the compactness of $V_{0}'$, we can find finitely many points $x_{0}^{0},\dots,x_{0}^{m}$ in $V_{0}'$ such that $V_{0}'\subseteq\bigcup\limits_{k=0}^{m}U_{x_{0}^{k}}$ with $\overline{U_{x_{0}^{k}}}\subseteq V_{0}$. Let $\mathcal{U}_{0}=\bigcup\limits_{k=0}^{m}U_{x_{0}^{k}}$, then $V_{0}'\subseteq\mathcal{U}_{0}$ and $\overline{\mathcal{U}_{0}}=\overline{\bigcup\limits_{k=0}^{m}U_{x_{0}^{k}}}=\bigcup\limits_{k=0}^{m}\overline{U_{x_{0}^{k}}}\subseteq V_{0}$.
 
Another step is to construct subsets $V_{i}'$ and $\mathcal{U}_{i}$ of $\mathcal{G}^{(0)}$ such that $U_{y}\subseteq V_{i}'$ when $y\in\partial V_{0}\cap V_{i}$ for some $i\in\{1,\dots,n\}$, and such that $V_{i}'\subseteq\mathcal{U}_{i}\subseteq\overline{\mathcal{U}_{i}}\subseteq V_{i}$. From the second step, we see that $U_{y}\subseteq V_{i}$ and $U_{y}\cap\widetilde{U_{z}}=\varnothing$, $\forall z\in\partial V_{i}$. For fixed $i$ and $z\in\partial V_{i}$, $z\notin V_{i}$ implies that there is at least one $l\in\{0,\dots,i-1,i+1,\dots,n\}$ such that $z\in V_{l}$. Therefore, there exists an open neighborhood $U_{z}'$ of $z$ such that $U_{z}'\subseteq V_{l}$ as $V_{l}$ is open. Similar to the second step, for any $h\in\partial V_{l}$, there exist open neighborhoods $\widetilde{U_{h}}$ of $h$ and $U_{z}''$ of $z$ such that $\widetilde{U_{h}}\cap U_{z}''=\varnothing$. Replacing $U_{z}'$ by $U_{z}'\cap U_{z}''$, we have $\widetilde{U_{h}}\cap U_{z}'=\varnothing$ and $U_{z}'\subseteq V_{l}$. Set $U_{z}=\widetilde{U_{z}}\cap U_{z}'$, then $U_{z}\subseteq V_{l}$, $U_{z}\cap U_{y}=\varnothing$ and $U_{z}\cap\widetilde{U_{h}}=\varnothing$. Let $V_{i}'=\overline{V_{i}}\backslash\bigcup\limits_{z\in\partial V_{i}}U_{z}$. Since $U_{z}\cap U_{y}=\varnothing$, we get the desired inclusion $U_{y}\subseteq V_{i}'$. Repeating the second step, we may obtain an open subset $\mathcal{U}_{i}$ such that $V_{i}'\subseteq\mathcal{U}_{i}\subseteq\overline{\mathcal{U}_{i}}\subseteq V_{i}$. 

Finally, we show that there are open subsets $\mathcal{U}_{0},\dots, \mathcal{U}_{n}$ of $\mathcal{G}^{(0)}$ covering $\mathcal{G}^{(0)}$ such that $\overline{\mathcal{U}_{i}}\subseteq V_{i}$ for each $i$. By the third step, for any $y\in\partial V_{0}$, there exist $i\in\{1,\dots,n\}$ and a closed subset $V_{i}'$ such that $U_{y}\subseteq V_{i}'$. Therefore, $\bigcup\limits_{y\in\partial V_{0}}U_{y}\subseteq\bigcup\limits_{i=1}^{n}V_{i}'\subseteq\bigcup\limits_{i=1}^{n}\mathcal{U}_{i}$. It follows that $\overline{V_{0}}=V_{0}'\cup\bigcup\limits_{y\in\partial V_{0}}U_{y}\subseteq\mathcal{U}_{0}\cup\bigcup\limits_{i=1}^{n}\mathcal{U}_{i}=\bigcup\limits_{i=0}^{n}\mathcal{U}_{i}$. Repeating the previous steps, we obtain 
$$
\overline{V_{j}}\subseteq\bigcup\limits_{i=0}^{n}\mathcal{U}_{i}, \text{ for each } j\in\{0,\dots,n\}.
$$
Hence, $\mathcal{G}^{(0)}\subseteq\bigcup\limits_{j=0}^{n}\overline{V_{j}}\subseteq\bigcup\limits_{i=0}^{n}\mathcal{U}_{i}$ with $\overline{\mathcal{U}_{i}}\subseteq V_{i}$ for each $i$.
\end{proof}

The following lemma constructs a nested family of open precompact subsets of $\mathcal{G}^{(0)}$ that satisfy several favorable properties. This is the main technical ingredient for the proof of the main proposition.

\begin{lemma}\label{lemma 5.6}
Let $\mathcal{G}$ be an \'etale groupoid with compact unit space, and with dynamic asymptotic dimension function $f$. Let $N:\mathbb{R}^{+}\rightarrow\mathbb{N}$ be a nondecreasing function. For any $R>0$, for the open precompact subset $K={\{z\in\mathcal{G}\mid\ell(z)< R\}}$ of $\mathcal{G}$, and each $i\in\{0,\dots,f(N(R)R+1)\}$, there is a nested family
$$
U^{(0)}_{i}\subseteq U^{(1)}_{i}\subseteq\cdots\subseteq U_{i}^{(N(R))}
$$
of open precompact subsets of $\mathcal{G}^{(0)}$ with the following properties.
\begin{enumerate}
\item The family $\{U_{0}^{(0)},\dots, U_{f(N(R)R+1)}^{(0)}\}$ covers $\mathcal{G}^{(0)}$;
\item For all $i,n$, $\overline{U^{(n)}_{i}}\subseteq U_{i}^{(n+1)}$;
\item For all $i,n$, 
$$
s(K\cap r^{-1}(U^{(n)}_{i}))\subseteq U_{i}^{(n+1)};
$$
\item For all $i$, the set
$$
\{g\in K\mid s(g),r(g)\in U_{i}^{(N(R))}\}
$$
generates an open precompact subgroupoid of $\mathcal{G}$.
\end{enumerate}
\end{lemma}
\begin{proof}
We prove this lemma in two steps. First, we build a nested sequence of open subsets of $\mathcal{G}^{(0)}$. Given $R>0$ and a nondecreasing function $N:\mathbb{R}^{+}\rightarrow\mathbb{N}$. Observe that $K$ is symmetric and contains $\mathcal{G}^{(0)}$. For each $n\in\mathbb{N}$, we define
\begin{equation*}
\overline{K}^{n}=\{g_{n}\cdots g_{1}\mid g_{k}\in\overline{K} \text{ and } s(g_{k+1})=r(g_{k}), \forall k\in\{1,\dots,n-1\}\}.
\end{equation*}
Since $\overline{K}$ is compact, it follows from Lemma \ref{lemma 4.11} that $\overline{K}^{n}$ is also compact. Furthermore, we have $\overline{K}^{n}\subseteq\overline{K}^{n+1}$ for each $n\in\mathbb{N}$, provided we set $g_{n+1}=r(g_{n})\in \overline{K}$. For any $z\in\overline{K}^{N(R)}$, we can observe that $\ell(z)\leq N(R)R<N(R)R+1$, thus $\overline{K}^{N(R)}$ is contained in $K'=\{z\in\mathcal{G}\mid\ell(z)<N(R)R+1\}$, which is an open precompact symmetric subset of $\mathcal{G}$ containing $\mathcal{G}^{(0)}$. According to Lemma \ref{lemma 5.4}, there exist open precompact subsets $V_{0},\dots, V_{f(N(R)R+1)}$ of $\mathcal{G}^{(0)}$ that cover $\mathcal{G}^{(0)}$ such that for each $i$, the set
$\{g\in K'\mid s(g),r(g)\in V_{i}\}$ generates an open precompact subgroupoid of $\mathcal{G}$.  

In addition, for each $i\in\{0,\dots,f(N(R)R+1)\}$ and each $n\in\{0,\dots,N(R)\}$, we define
$$
V_{i}^{(n)}:=\{x\in V_{i}\mid s(\overline{K}^{N(R)-n}\cap r^{-1}(x))\subseteq V_{i}\}.
$$
The sequence $\mathcal{G}^{(0)}=\overline{K}^{0}\subseteq\overline{K}\subseteq\cdots\subseteq\overline{K}^{N(R)}$ induces a nested sequence of subsets
$$
V^{(0)}_{i}\subseteq V_{i}^{(1)}\subseteq\cdots\subseteq V_{i}^{(N(R))}=V_{i}.
$$
By Lemma \ref{lemma 5.4}, for any $x\in\mathcal{G}^{(0)}$, $s(K'\cap r^{-1}(x))\subseteq V_{i}$ for some $i$, namely, if $g\in K'$ with $r(g)=x$, then $s(g)\in V_{i}$. By Lemma \ref{lemma 5.3}, $x=r(g)\in V_{i}$. Further, $s(\overline{K}^{N(R)}\cap r^{-1}(x))\subseteq s(K'\cap r^{-1}(x))\subseteq V_{i}$, then $x\in V_{i}^{(0)}$. Therefore, the family $\{V_{0}^{(0)},\dots, V^{(0)}_{f(N(R)R+1)}\}$ covers $\mathcal{G}^{(0)}$. 

We claim that $V_{i}^{(n)}$ is open. Say for contradiction that $V_{i}^{(n)}$ is not open, then we may find a net $(x_{\lambda})_{\lambda\in\Lambda}$ in $\mathcal{G}^{(0)}$ with $x_{\lambda}\rightarrow x$ for some $x\in V_{i}^{(n)}$ such that for each $\lambda\in\Lambda$, there exists an element $g_{\lambda}\in \overline{K}^{N(R)-n}\cap r^{-1}(x_{\lambda})$ such that $s(g_{\lambda})\notin V_{i}$. Since $\overline{K}^{N(R)-n}$ is compact, the net $(g_{\lambda})_{\lambda\in\Lambda}$ has a convergent subnet $(g_{\lambda_{j}})_{j\in J}$, thus there exists $g\in\overline{K}^{N(R)-n}$ such that $\lim\limits_{j}g_{\lambda_{j}}=g$. As $r$ is continuous and $r(g_{\lambda})=x_{\lambda}$, $x=\lim\limits_{\lambda}x_{\lambda}=\lim\limits_{j}x_{\lambda_{j}}=\lim\limits_{j}r(g_{\lambda_{j}})=r(g)$, then $g\in\overline{K}^{N(R)-n}\cap r^{-1}(x)$. Because $V_{i}$ is open and $s(g_{\lambda})\notin V_{i}$, $s(g)\notin V_{i}$, which contradicts that $x\in V_{i}^{(n)}$. Hence, $V_{i}^{(n)}$ is open as claimed.

In the second step, we construct the desired nested sequence $U_{i}^{(n)}$. Since $\mathcal{G}^{(0)}$ is compact, by Lemma \ref{lemma 5.5}, there exist open subsets $U_{0}^{(0)},\dots, U_{f(N(R)R+1)}^{(0)}$ of $\mathcal{G}^{(0)}$ covering $\mathcal{G}^{(0)}$ and 
$\overline{U_{i}^{(0)}}\subseteq V_{i}^{(0)}$ for each $i$. For each $i$, 
$$
\overline{s(K\cap r^{-1}(U_{i}^{(0)}))}\subseteq s(\overline{K}\cap r^{-1}(\overline{U_{i}^{(0)}}))\subseteq s(\overline{K}\cap r^{-1}(V_{i}^{(0)}))\subseteq V_{i}^{(1)}.
$$ 
By normality of the compact set $\overline{V_{i}}$, there exists an open set $U_{i}^{(1)}$ such that 
$$
\overline{s(K\cap r^{-1}(U_{i}^{(0)}))}\subseteq U_{i}^{(1)}\subseteq\overline{U_{i}^{(1)}}\subseteq V_{i}^{(1)}.
$$
We also see that
$$
\overline{s(K\cap r^{-1}(U_{i}^{(1)}))}\subseteq s(\overline{K}\cap r^{-1}(V_{i}^{(1)}))\subseteq V_{i}^{(2)},
$$
there exists an open set $U_{i}^{(2)}$ such that
$$
\overline{s(K\cap r^{-1}(U_{i}^{(1)}))}\subseteq U_{i}^{(2)}\subseteq\overline{U_{i}^{(2)}}\subseteq V_{i}^{(2)}.
$$
Repeating this process, we get a nested sequence of open subsets
$$
U_{i}^{(0)}\subseteq U_{i}^{(1)}\subseteq\cdots\subseteq U_{i}^{(N(R))},
$$
such that for each $i,n$,
$$
\overline{s(K\cap r^{-1}(U_{i}^{(n)}))}\subseteq U_{i}^{(n+1)}\text{ and }\overline{U_{i}^{(n)}}\subseteq V_{i}^{(n)}.
$$
Since $U_{i}^{(n)}\subseteq V_{i}^{(n)}\subseteq V_{i}\subseteq\mathcal{G}^{(0)}$, the set $U_{i}^{(n)}$ is precompact. Moreover, 
\begin{equation*}
\overline{U_{i}^{(n)}}\subseteq\overline{s(K\cap r^{-1}(U_{i}^{(n)}))}\subseteq U_{i}^{(n+1)}.
\end{equation*}
Note that for all $i$, the set
\begin{equation*}
\{g\in K\mid s(g),r(g)\in U_{i}^{(N(R))}\}\subseteq\{g\in K'\mid s(g),r(g)\in V_{i}\},
\end{equation*}
which implies that the set $\{g\in K\mid s(g),r(g)\in U_{i}^{(N(R))}\}$ generates an open precompact subgroupoid of $\mathcal{G}$.
\end{proof}

We have established the necessary lemmas for Proposition \ref{prop 5.1}, enabling us to commence its proof. 
 
\begin{proof}[Proof of Proposition \ref{prop 5.1}] 
\textbf{Step 1}. In the first step, we build a desired open cover of $\mathcal{G}^{(0)}$. Given $R>0$, the open precompact subset $K=\{z\in\mathcal{G}\mid \ell(z)<R\}$ of $\mathcal{G}$, and $\varepsilon>0$. Clearly, $K=K^{-1}$, $s(K)\cup r(K)=\mathcal{G}^{(0)}\subseteq K$. Let $p=[\frac{\alpha}{1-\alpha}]+1$, and let $N(R)=[(R+1)^{c}]$ for some constant $c$ in Lemma \ref{lemma 5.6}, such that
\begin{equation}\label{eq 5.1}
\frac{2p\left(f(N(R)R+1)+1\right)+2p\left(f(N(R)R+1)+1\right)^{\frac{1}{p}+1}}{N(R)}<\varepsilon,
\end{equation}
such constant $c$ exists and depends on given parameters $\alpha, R,\varepsilon$. Let $\{U_{i}^{(n)}\mid i\in\{0,\dots,f(N(R)R+1)\}, n\in\{0,\dots, N(R)\}\}$ have the properties in Lemma \ref{lemma 5.6}. For each $i$, let $U_{i}=U_{i}^{(N(R))}$. By Lemma \ref{lemma 5.6} (1)(2), $\{U_{0}^{(0)},\dots, U_{f(N(R)R+1)}^{(0)}\}$ covers $\mathcal{G}^{(0)}$ and $U_{i}^{(0)}\subseteq U_{i}^{(N(R))}=U_{i}$, thus the family $\{U_{0},\dots, U_{f(N(R)R+1)}\}$ covers $\mathcal{G}^{(0)}$. By Lemma \ref{lemma 5.6} (3), we have that $s(K\cap r^{-1}(U_{i}^{(0)}))\subseteq U_{i}^{(1)}\subseteq U_{i}$ for some $i$, which implies that for any $x\in\mathcal{G}^{(0)}$, $s(K\cap r^{-1}(x))\subseteq U_{i}$ for some $U_{i}$. By Lemma \ref{lemma 5.6} (4), we see that
$\{g\in K\mid s(g), r(g)\in U_{i}\} $ generates an open precompact subgroupoid of $\mathcal{G}$. 

\textbf{Step 2}. The second step is to construct a partition of unity. By Lemma \ref{lemma 5.6} (2), for each $i\in\{0,\dots,f(N(R)R+1)\}$, and each $n\in\{1,\dots, N(R)\}$, $\overline{U_{i}^{(n-1)}}\subseteq U_{i}^{(n)}$. According to Urysohn's lemma, there exist continuous functions
$$
\psi_{i}^{(n)}:\mathcal{G}^{(0)}\rightarrow [0,1]
$$
such that $\psi_{i}^{(n)}(x)=1$, for all $x\in U_{i}^{(n-1)}$; $\psi^{(n)}_{i}(x)=0$, for all $ x\in\mathcal{G}^{(0)}\backslash U_{i}^{(n)}$. 

For $i\in\{0,\dots,f(N(R)R+1)\}$, put
$$
\psi_{i}=\frac{1}{N(R)}\sum\limits_{n=1}^{N(R)}\psi_{i}^{(n)},
$$
then $\supp(\psi_{i})\subseteq U_{i}$ and $\psi_{i}(x)\leq 1, \text{ for all }x\in\mathcal{G}^{(0)}$. 
Define
$$
\phi_{i}=\frac{\psi_{i}}{\left(\sum\limits_{j=0}^{f(N(R)R+1)}\psi_{j}^{p}\right)^{\frac{1}{p}}},
$$ 
which are continuous functions and $\supp(\phi_{i})\subseteq U_{i}$.
By Lemma \ref{lemma 5.6} (1), for any $x\in\mathcal{G}^{(0)}$, there exists at least one $j\in\{0,\dots,f(N(R)R+1)\}$ such that $x\in U_{j}^{(0)}$, thus $\psi_{j}(x)=1$, so
$\left(\sum\limits_{j=0}^{f(N(R)R+1)}\psi_{j}^{p}(x)\right)^{\frac{1}{p}}\geq 1$, therefore $$0\leq\phi_{i}(x)\leq 1 \text{ and } 
\sum\limits_{i=0}^{f(N(R)R+1)}\left(\phi_{i}(x)\right)^{p}=1. $$

\textbf{Step 3}. The final step is to show that for any $g\in K$, 
$$p\cdot\sum\limits_{i=0}^{f(N(R)R+1)}\left|\phi_{i}(s(g))-\phi_{i}(r(g))\right|<\varepsilon. $$
For each $j\in\{0,\dots,f(N(R)R+1)\}$, set $U_{j}^{(n)}=\mathcal{G}^{(0)}$ for $n\geq N(R)+1$.
Define $M=M_{j}:=\min\{n\mid r(g)\in U_{j}^{(n)}\}$. Thus $r(g)\in U_{j}^{(M)}\backslash U_{j}^{(M-1)}$, and
$$
\psi_{j}^{(n)}(r(g))=\begin{cases} 1,\quad & n\in [M+1, N(R)]\\ 0,\quad & n\in [1, M-1]\\ [0,1],\quad & n=M.\end{cases}
$$
Hence
\begin{equation}\label{eq 5.2}
\frac{N(R)-(M+1)}{N(R)}\leq\psi_{j}(r(g))\leq\frac{N(R)-M}{N(R)}.
\end{equation}
Since $g\in K$ and $r(g)\in U_{j}^{(M)}$, $g\in K\cap r^{-1}(U_{j}^{(M)})$. By Lemma \ref{lemma 5.6} (2), we see that $s(g)\in s(K\cap r^{-1}(U_{j}^{(M)}))\subseteq U_{j}^{(M+1)}$. Note that $r(g)\notin U_{j}^{(M-1)}\supseteq s(K\cap r^{-1}(U_{j}^{(M-2)}))$, thus $s(g^{-1})=r(g)\notin s(K\cap r^{-1}(U_{j}^{(M-2)}))$. We see that $g^{-1}\in K$ as $K=K^{-1}$, thus $g^{-1}\notin r^{-1}(U_{j}^{(M-2)})$, which implies that $s(g)=r(g^{-1})\notin U_{j}^{(M-2)}$. Therefore $s(g)\in U_{j}^{(M+1)}\backslash U_{j}^{(M-2)}$, and
$$
\psi_{j}^{(n)}(s(g))=\begin{cases} 1,\quad & n\in [M+2, N(R)]\\ 0,\quad & n\in [1, M-2]\\ [0,1],\quad & n=M-1, M, M+1.\end{cases}
$$
Hence
\begin{equation}\label{eq 5.3}
\frac{N(R)-(M+2)}{N(R)}\leq\psi_{j}(s(g))\leq\frac{N(R)-(M-1)}{N(R)}.
\end{equation}
Combining equation \eqref{eq 5.2} and equation \eqref{eq 5.3}, we have that
\begin{equation}\label{eq 5.4}
\left|\psi_{j}(s(g))-\psi_{j}(r(g))\right|\leq\frac{2}{N(R)}. 
\end{equation}
By Lemma \ref{lemma 5.3}, there exists at least one $j\in\{0,\dots,f(N(R)R+1)\}$ such that $s(g), r(g)$ in $U_{j}^{(0)}$, hence $\psi_{j}(r(g))=\psi_{j}(s(g))=1$, which implies that
\begin{equation}\label{eq 5.5}
\sum\limits_{j=0}^{f(N(R)R+1)}\psi_{j}^{p}(r(g))\geq 1 \text{ and } \sum\limits_{j=0}^{f(N(R)R+1)}\psi_{j}^{p}(s(g))\geq 1.
\end{equation}
For fixed $i\in\{0,\dots,f(N(R)R+1)\}$ and $g\in K$,
\begin{align}\label{eq 5.6}
\left|\phi_{i}(r(g))-\phi_{i}(s(g))\right|&=\left|\frac{\psi_{i}(r(g))}{\left(\sum\limits_{j=0}^{f(N(R)R+1)}\psi_{j}^{p}(r(g))\right)^{\frac{1}{p}}}-\frac{\psi_{i}(s(g))}{\left(\sum\limits_{j=0}^{f(N(R)R+1)}\psi_{j}^{p}(s(g))\right)^{\frac{1}{p}}}\right|\nonumber\\
&\leq\frac{1}{\left(\sum\limits_{j=0}^{f(N(R)R+1)}\psi_{j}^{p}(r(g))\right)^{\frac{1}{p}}}\left|\psi_{i}(r(g))-\psi_{i}(s(g))\right|\nonumber
\\&+\left|\psi_{i}(s(g))\right|\left|\frac{1}{\left(\sum\limits_{j=0}^{f(N(R)R+1)}\psi_{j}^{p}(r(g))\right)^{\frac{1}{p}}}-\frac{1}{\left(\sum\limits_{j=0}^{f(N(R)R+1)}\psi_{j}^{p}(s(g))\right)^{\frac{1}{p}}}\right|.
\end{align}
Combining equation \eqref{eq 5.4} and equation \eqref{eq 5.5}, we have
\begin{equation}\label{eq 5.7}
\frac{1}{\left(\sum\limits_{j=0}^{f(N(R)R+1)}\psi_{j}^{p}(r(g))\right)^{\frac{1}{p}}}\left|\psi_{i}(r(g))-\psi_{i}(s(g))\right|\leq\frac{2}{N(R)}.
\end{equation}
On the other hand, using equations \eqref{eq 5.4} and \eqref{eq 5.5} again, we obtain
\begin{align}\label{eq 5.8}
&\left|\psi_{i}(s(g))\right|\left|\frac{1}{\left(\sum\limits_{j=0}^{f(N(R)R+1)}\psi_{j}^{p}(r(g))\right)^{\frac{1}{p}}}-\frac{1}{\left(\sum\limits_{j=0}^{f(N(R)R+1)}\psi_{j}^{p}(s(g))\right)^{\frac{1}{p}}}\right|\nonumber\\
&\leq \frac{\left|\left(\sum\limits_{j=0}^{f(N(R)R+1)}\psi_{j}^{p}(s(g))\right)^\frac{1}{p}-\left(\sum\limits_{j=0}^{f(N(R)R+1)}\psi_{j}^{p}(r(g))\right)^{\frac{1}{p}}\right|}{\left(\sum\limits_{j=0}^{f(N(R)R+1)}\psi_{j}^{p}(r(g))\right)^{\frac{1}{p}}\left(\sum\limits_{j=0}^{f(N(R)R+1)}\psi_{j}^{p}(s(g))\right)^{\frac{1}{p}}}
\nonumber\\
&\leq\left(\sum\limits_{j=0}^{f(N(R)R+1)}\left|\psi_{j}(s(g))-\psi_{j}(r(g))\right|^{p}\right)^{\frac{1}{p}}\nonumber\\
&\leq\frac{2\left(f(N(R)R+1)+1\right)^{\frac{1}{p}}}{N(R)},
\end{align}
where the second inequality comes from the Minkowski inequality $$\left|\left(\sum\limits_{j=0}^{d}a_{j}^{p}\right)^\frac{1}{p}-\left(\sum\limits_{j=0}^{d}b_{j}^{p}\right)^{\frac{1}{p}}\right|\leq\left(\sum\limits_{j=0}^{d}|a_{j}-b_{j}|^{p}\right)^\frac{1}{p}, \text{ for }a_{j}, b_{j}\geq 0.$$
According to equations \eqref{eq 5.6}, \eqref{eq 5.7}, and \eqref{eq 5.8}, we have
\begin{align*}
\left|\phi_{i}(r(g))-\phi_{i}(s(g))\right|&\leq\frac{2}{N(R)}+\frac{2\left(f(N(R)R+1)+1\right)^{\frac{1}{p}}}{N(R)}\\
&=\frac{2+2\left(f(N(R)R+1)+1\right)^{\frac{1}{p}}}{N(R)}.
\end{align*}
Finally, from the equation \eqref{eq 5.1}, we obtain
\begin{align*}
&p\cdot\sum\limits_{i=0}^{f(N(R)R+1)}\left|\phi_{i}(r(g))-\phi_{i}(s(g))\right|\\
&\leq\frac{2p\left(f(N(R)R+1)+1\right)+2p\left(f(N(R)R+1)+1\right)^{\frac{1}{p}+1}}{N(R)}
<\varepsilon. 
\end{align*}
\end{proof}

\begin{remark}\label{rmk 5.7}
Let $0<\alpha<1$, $R>0$ and $\varepsilon>0$. We consider the function $f(x)=[x^{\alpha}]$, $p=[\frac{\alpha}{1-\alpha}]+1$ to examine the existence of the constant $c$. Through precise calculations, we found that if the constant $c$ satisfies the following inequality:
$$
c>\max\left\{\frac{\alpha(p+1)\In 2-p\In \frac{\varepsilon}{4p}}{\left(p-\alpha(p+1)\right)\In (R+1)},\frac{\alpha\In 2-\In\frac{\varepsilon}{8p}}{(1-\alpha)\In (R+1)}\right\},
$$
we have 
$$\frac{2p\left(f([(R+1)^{c}]R+1)+1\right)+2p\left(f([(R+1)^{c}]R+1)+1\right)^{\frac{1}{p}+1}}{[(R+1)^{c}]}<\varepsilon.$$
\end{remark}

After thoroughly examining the proof of Proposition 7.1 in \cite{GWY17}, we conclude that finite dynamic asymptotic dimension leads to an almost invariant partition of unity without the restriction of the compactness of the unit space. Consequently, we present a stronger result as follows:

\begin{proposition}[\cite{GWY17}]\label{prop 5.8}
Let $\mathcal{G}$ be an \'etale groupoid with dynamic asymptotic dimension $d$. Then for any open precompact subset $K$ of $\mathcal{G}$ and $\varepsilon>0$, there exists an open cover $\{U_{0},\dots, U_{d}\}$ of $r(K)\cup s(K)$ with the following properties.
\begin{enumerate}
\item For each $i$, the set $$\{g\in K\mid s(g), r(g)\in U_{i}\}$$ is contained in an ( open and ) precompact subgroupoid of $\mathcal{G}$.
\item For each $x\in r(K)\cup s(K)$, the `partial orbit' $s(r^{-1}(x)\cap K)$ is completely contained in some $U_{i}$.
\item There exists a collection of continuous functions $\{\phi_{i}:\mathcal{G}^{(0)}\rightarrow [0,1]\}_{i=0}^{d}$ such that $\supp(\phi_{i})\subseteq U_{i}$ for each $i$, $\sum\limits_{i=0}^{d}\phi_{i}(x)\leq 1$ for all $x\in\mathcal{G}^{(0)}$, and $\sum\limits_{i=0}^{d}\phi_{i}(x)=1$ for all $x\in r(K)\cup s(K)$, and such that for any $g\in K$ and each $i$,
$$
\vert\phi_{i}(r(g))-\phi_{i}(s(g))\vert<\varepsilon.
$$
\end{enumerate}
\end{proposition}

\section{An application to amenability}

Our motivation for investigating the relationship between the dynamic asymptotic dimension growth and the amenability of groupoids stems from three key prior results. First, Guentner, Willett, and Yu established that the finite dynamic asymptotic dimension of a free, locally compact Hausdorff \'etale groupoid implies amenability \cite{GWY17}. Second, they demonstrated that the finite dynamical complexity of a locally compact Hausdorff \'etale groupoid also implies amenability \cite{GWY24}. It is worth noting that finite dynamic asymptotic dimension implies finite dynamical complexity when the former is either zero or one. Finally, Corollary \ref{cor 4.21} from Section \ref{section: dad} shows that for a discrete metric space $X$ of bounded geometry, the coarse groupoid $G(X)$ is amenable if it has subexponential dynamic asymptotic dimension growth. 

Building on these foundational insights, this section establishes broader sufficient conditions for the amenability of groupoids. We show that every $\sigma$-compact, locally compact, Hausdorff, and \'etale groupoid $\mathcal{G}$ with compact unit space and dynamic asymptotic dimension growth at most $x^{\alpha}$ $(0<\alpha<1)$ is amenable (Theorem \ref{th 6.4}), thus answering Question \ref{q 1.2} posed in the introduction. A direct consequence is that the Baum--Connes conjecture holds for such groupoids (Corollary \ref{cor 6.6}). Furthermore, applying Corollary \ref{cor 6.6} to [\cite{PY26}, Corollary 5.6] provides a sufficient condition for the rational HK conjecture (Corollary \ref{cor 6.10}). Finally, the section concludes with a discussion of several open problems.

The following definition from [\cite{GWY24}, Definition A.8] is equivalent to a slight variation of Definition 5.6.13 in \cite{BO08} (specifically addressing the \'etale case) and provides the formulation used in the proofs of the two theorems presented in this section. For a standard treatment of amenable groupoids, we refer the reader to \cite{AR00}.
\begin{definition}[\cite{GWY24}]
Let $\mathcal{G}$ be a locally compact, Hausdorff, \'etale groupoid. Say $\mathcal{G}$ is amenable if for all compact subset $K\subseteq\mathcal{G}$ and all $\varepsilon>0$, there exists a continuous, compactly supported function $\mu:\mathcal{G}\rightarrow [0,1]$ such that:
	\begin{enumerate}
\item for all $x\in\mathcal{G}^{(0)}$, we have $\sum\limits_{g\in\mathcal{G}_{x}}\mu(g)\leq 1$; 
\item for all $k\in K$, we have $\left| 1-\sum\limits_{g\in\mathcal{G}_{r(k)}}\mu(g)\right|<\varepsilon$;
\item for all $k\in K$, we have $\sum\limits_{g\in\mathcal{G}_{r(k)}}\left|\mu(g)-\mu(gk)\right|<\varepsilon$.
	\end{enumerate}
\end{definition}

The following lemma demonstrates that any open precompact subgroupoid of an \'etale groupoid is amenable. This property plays a crucial role in the proofs of the two theorems in this section.
\begin{lemma}[\cite{GWY24}]\label{lemma 6.2}
Let $\mathcal{G}$ be an \'etale groupoid, and let $H$ be an open subgroupoid of $\mathcal{G}$ with compact closure. Then $H$ is amenable.
\end{lemma}

The following theorem is inspired by the work of Guentner, Willett, and Yu, particularly their proof that groupoids with finite dynamical complexity are amenable \cite{GWY24}. In contrast to their approach in \cite{GWY17}, which relied on the nuclear dimension of the reduced groupoid $C^{*}$-algebra under a freeness assumption, we employ a direct method from [\cite{GWY24}, Theorem A.9] in conjunction with Proposition \ref{prop 5.8}. This allows us to strengthen the conclusion of [\cite{GWY17}, Corollary 8.25] by eliminating the freeness assumption.

\begin{theorem}\label{th 6.3}
Let $\mathcal{G}$ be a locally compact, Hausdorff, and \'etale groupoid with finite dynamic asymptotic dimension. Then $\mathcal{G}$ is amenable.
\end{theorem}

\begin{proof}
Let us assume that $\mathcal{G}$ has dynamic asymptotic dimension $d$. For any compact subset $K$ of $\mathcal{G}$ and any $\varepsilon>0$. Using local compactness, we expand $K$ slightly and write it as $K'$, which is open and precompact. Replacing $K'$ with $K'\cup K'^{-1}\cup r(K')\cup s(K')$, we may assume that $K'=K'\cup K'^{-1}\cup r(K')\cup s(K')$, which remains open and precompact. Next, let $\{U_{0},\dots,U_{d}\}$ be an open cover of $r(K')\cup s(K')$, satisfying the properties outlined in Proposition \ref{prop 5.8} for the open precompact subset $K'$ and the error estimate $\frac{\varepsilon}{d+2}$. By taking the intersection of each $U_{i}$ with $s(K')\cup r(K')$, we may assume that each $U_{i}$ is contained in $s(K')\cup r(K')\subseteq K'$. Let $\mathcal{G}_{i}$ represent the subgroupoid generated by the set $K_{i}=\{k\in K'\mid s(k),r(k)\in U_{i}\}$. According to Proposition \ref{prop 5.8}, $\mathcal{G}_{i}$ is both open and precompact. 

We claim that $K'=\bigcup\limits_{i}K_{i}$, where each $K_{i}$ is open, symmetric, precompact, and has compact closure in $\mathcal{G}_{i}$. We first demonstrate that $K'=\bigcup\limits_{i}K_{i}$. According to items (ii) and (iii) of [\cite{GWY17}, Lemma 7.4], we have $s(K'\cap r^{-1}(U_{i}^{(0)}))\subseteq U_{i}^{(1)}$ and $U_{i}^{(0)}\subseteq U_{i}^{(1)}\subseteq U_{i}$. Thus, it follows that $K'\cap r^{-1}(U_{i}^{(0)})\subseteq K_{i}$. Furthermore, by the item  (i) of [\cite{GWY17}, Lemma 7.4], the collection $\{U_{0}^{(0)},\dots, U_{d}^{(0)}\}$ covers $\overline{r(K')\cup s(K')}$. Therefore, 
\begin{align*}
\bigcup_{i=0}^{d}K_{i}\supseteq \bigcup_{i=0}^{d}K'\cap r^{-1}(U_{i}^{(0)})&=K'\cap r^{-1}\left(\bigcup_{i=0}^{d}U_{i}^{(0)}\right)\\&\supseteq K'\cap r^{-1}\left(\overline{r(K')\cup s(K')}\right)\supseteq K'.
\end{align*}
Clearly, $K_{i}$ is open, precompact, and has compact closure in $\mathcal{G}_{i}$. It remains to show that $K_{i}=K_{i}^{-1}$. Indeed, for any $k\in K_{i}$, specifically, $k\in K'$, $s(k),r(k)\in U_{i}$. Consequently, $k^{-1}\in K'$ since $K'=K'^{-1}$. Furthermore, since $r(k^{-1})=s(k)$ and $s(k^{-1})=r(k)$, we have $r(k^{-1}), s(k^{-1})\in U_{i}$, which implies that $k^{-1}\in K_{i}$. This completes the proof of the claim.

Since $\mathcal{G}_{i}$ is an open precompact subgroupoid of $\mathcal{G}$, it follows from Lemma \ref{lemma 6.2} that $\mathcal{G}_{i}$ is amenable. For each $i$, let $\mu_{i}:\mathcal{G}_{i}\rightarrow [0,1]$ be a continuous, compactly supported function as defined in the context of amenability concerning the compact subset $\overline{K_{i}}$ and the error estimate $\frac{\varepsilon}{d+2}$. We extend $\mu_{i}$ to the entire $\mathcal{G}$ by defining
$\mu_{i}(x)=0$ when $x\notin\mathcal{G}_{i}$. 
We then define a function
$\mu:\mathcal{G}\rightarrow [0,1]$,
$g\mapsto\sum\limits_{i=0}^{d}\phi_{i}(s(g))\mu_{i}(g)$. Clearly, $\mu$ is a well-defined function that is continuous and has compact support in $\mathcal{G}$.

For any $x\in\mathcal{G}^{(0)}$, 
\begin{align*}
\sum\limits_{g\in\mathcal{G}_{x}}\mu(g)&=\sum\limits_{g\in\mathcal{G}_{x}}\sum\limits_{i=0}^{d}\phi_{i}(s(g))\mu_{i}(g)\\
&=\sum\limits_{i=0}^{d}\phi_{i}(x)\sum\limits_{g\in(\mathcal{G}_{i})_{x}}\mu_{i}(g)\\
&\leq\sum_{i=0}^{d}\phi_{i}(x)\leq 1, 
\end{align*}
where the first inequality comes from the amenability of $\mathcal{G}_{i}$. The final inequality is due to the item (3) of Proposition \ref{prop 5.8}.

For any $k\in K\subseteq K'$, 
\begin{align*}
\left|1-\sum_{g\in\mathcal{G}_{r(k)}}\mu(g)\right|&=\left|1-\sum_{g\in\mathcal{G}_{r(k)}}\sum_{i=0}^{d}\phi_{i}(s(g))\mu_{i}(g)\right|\\
&=\left|\sum_{i=0}^{d}\phi_{i}(r(k))\left(1-\sum_{g\in(\mathcal{G}_{i})_{r(k)}}\mu_{i}(g)\right)\right|\\
&\leq\sum_{i=0}^{d}\phi_{i}(r(k))\left|1-\sum_{g\in(\mathcal{G}_{i})_{r(k)}}\mu_{i}(g)\right|\\&<\sum_{i=0}^{d}\phi_{i}(r(k))\varepsilon=\varepsilon,
\end{align*}
where the strict inequality follows from the amenability of $\mathcal{G}_{i}$. The final equality is due to the partition of unity stated in Proposition \ref{prop 5.8}.

For any $k\in K\subseteq K'$, 
\begin{align*}
&\sum\limits_{g\in\mathcal{G}_{r(k)}}\left|\mu(g)-\mu(gk)\right|\\
&=\sum\limits_{g\in\mathcal{G}_{r(k)}}\left|\sum\limits_{i=0}^{d}\left(\phi_{i}(s(g))\mu_{i}(g)-\phi_{i}(s(gk))\mu_{i}(gk)\right)\right|\\
&\leq\sum\limits_{i=0}^{d}\sum\limits_{g\in\mathcal{G}_{r(k)}}\phi_{i}(r(k))\left|\mu_{i}(g)-\mu_{i}(gk)\right|+\mu_{i}(gk)\left|\phi_{i}(r(k))-\phi_{i}(s(k))\right|\\
&=\sum\limits_{i=0}^{d}\phi_{i}(r(k))\sum\limits_{g\in(\mathcal{G}_{i})_{r(k)}}\left|\mu_{i}(g)-\mu_{i}(gk)\right|\\&+\sum\limits_{i=0}^{d}\left|\phi_{i}(r(k))-\phi_{i}(s(k))\right|\sum\limits_{g\in(\mathcal{G}_{i})_{r(k)}}\mu_{i}(gk)\\&<\sum\limits_{i=0}^{d}\phi_{i}(r(k))\cdot\frac{\varepsilon}{d+2}+(d+1)\cdot\frac{\varepsilon}{d+2}\cdot 1=\varepsilon,
\end{align*}
where the strict inequality follows from the amenability of $\mathcal{G}_{i}$ and Proposition \ref{prop 5.8}.
\end{proof}

The following main theorem provides a definitive answer to Question \ref{q 1.2} by establishing that slow dynamic asymptotic dimension growth implies the amenability of a groupoid. Our proof adopts an approach analogous to that of Theorem \ref{th 6.3}, using the key technical Proposition \ref{prop 5.1}. As a generalization, the theorem extends Corollary \ref{cor 4.21} to general \'etale groupoids with slow growth.

\begin{theorem}\label{th 6.4}
Let $\mathcal{G}$ be a $\sigma$-compact, locally compact, Hausdorff, and \'etale groupoid with compact unit space and dynamic asymptotic dimension growth $f\preceq x^{\alpha} (0<\alpha<1)$. Then $\mathcal{G}$ is amenable.
\end{theorem}

\begin{proof}
Without loss of generality, we may assume that $\mathcal{G}$ has dynamic asymptotic dimension function $f\preceq x^{\alpha}$ for some $\alpha\in (0,1)$. For any compact subset $K$ of $\mathcal{G}$ and any $\varepsilon>0$, there exists $R>0$ such that $K\subseteq K'=\{z\in\mathcal{G}\mid\ell(z)<R\}$ as the length function $\ell$ is controlled. Note that $K'$ is open and precompact, $K'=K'^{-1}$, and $s(K')\cup r(K')=\mathcal{G}^{(0)}\subseteq K'$. By Proposition \ref{prop 5.1}, there exists a constant $c$ depending on parameters $\alpha, R, \varepsilon$; a positive integer $p$ depending only on $\alpha$; and an open cover $\{U_{0},\dots, U_{f([(R+1)^{c}]R+1)}\}$ of $\mathcal{G}^{(0)}$ with the properties for the open precompact subset $K'$ and the error estimate $\frac{\varepsilon}{2}$.  

Let $\mathcal{G}_{i}$ be the subgroupoid generated by the subset $K_{i}=\{k\in K'\mid s(k), r(k)\in U_{i}\}$. Similar to the claim in the proof of Theorem \ref{th 6.3}, we also have the following statement: $K'=\bigcup\limits_{i}K_{i}$ and $K_{i}$ is open, symmetric, precompact, and has compact closure in $\mathcal{G}_{i}$. According to Proposition \ref{prop 5.1}, $\mathcal{G}_{i}$ is an open precompact subgroupoid of $\mathcal{G}$. Therefore, by Lemma \ref{lemma 6.2}, $\mathcal{G}_{i}$ is amenable. For each $i$, let $\mu_{i}:\mathcal{G}_{i}\rightarrow [0,1]$ be a continuous, compactly supported function as in the definition of amenability for the compact subset $\overline{K_{i}}$ and the error estimate $\frac{\varepsilon}{2}$. We extend $\mu_{i}$ to the whole $\mathcal{G}$ by defining
$\mu_{i}(x)=0$, if $x\notin\mathcal{G}_{i}$. We then define
$\mu:\mathcal{G}\rightarrow [0,1], g\mapsto\sum\limits_{i=0}^{f([(R+1)^{c}]R+1)}\phi_{i}^{p}(s(g))\mu_{i}(g)$, which is a well-defined continuous, compactly supported function on $\mathcal{G}$. 

For any $x\in\mathcal{G}^{(0)}$, 
\begin{align*}
\sum\limits_{g\in\mathcal{G}_{x}}\mu(g)&=\sum\limits_{g\in\mathcal{G}_{x}}\sum\limits_{i=0}^{f([(R+1)^{c}]R+1)}\phi_{i}^{p}(s(g))\mu_{i}(g)\\
&=\sum\limits_{i=0}^{f([(R+1)^{c}]R+1)}\phi_{i}^{p}(x)\sum\limits_{g\in(\mathcal{G}_{i})_{x}}\mu_{i}(g)\\
&\leq\sum_{i=0}^{f([(R+1)^{c}]R+1)}\phi_{i}^{p}(x)=1,
\end{align*}
where the inequality follows from the amenability of $\mathcal{G}_{i}$, and the final equality is due to the partition of unity in Proposition \ref{prop 5.1}.

For any $k\in K\subseteq K'$, 
\begin{align*}
\left|1-\sum_{g\in\mathcal{G}_{r(k)}}\mu(g)\right|&=\left|1-\sum_{g\in\mathcal{G}_{r(k)}}\sum_{i=0}^{f([(R+1)^{c}]R+1)}\phi_{i}^{p}(s(g))\mu_{i}(g)\right|\\
&=\left|\sum_{i=0}^{f([(R+1)^{c}]R+1)}\phi_{i}^{p}(r(k))\left(1-\sum_{g\in(\mathcal{G}_{i})_{r(k)}}\mu_{i}(g)\right)\right|\\
&\leq\sum_{i=0}^{f([(R+1)^{c}]R+1)}\phi_{i}^{p}(r(k))\left|1-\sum_{g\in(\mathcal{G}_{i})_{r(k)}}\mu_{i}(g)\right|\\&<\sum_{i=0}^{f([(R+1)^{c}]R+1)}\phi_{i}^{p}(r(k))\cdot\frac{\varepsilon}{2}=\frac{\varepsilon}{2}<\varepsilon,
\end{align*}
where the second inequality follows from the amenability of $\mathcal{G}_{i}$, and the final equality comes from the partition of unity in Proposition \ref{prop 5.1}.

For any $k\in K\subseteq K'$, 
\begin{align*}
&\sum\limits_{g\in\mathcal{G}_{r(k)}}\left|\mu(g)-\mu(gk)\right|\\
&\leq\sum\limits_{g\in\mathcal{G}_{r(k)}}\left|\sum\limits_{i=0}^{f([(R+1)^{c}]R+1)}\left(\phi_{i}^{p}(s(g))\mu_{i}(g)-\phi_{i}^{p}(s(gk))\mu_{i}(gk)\right)\right|\\
&\leq\sum\limits_{i=0}^{f([(R+1)^{c}]R+1)}\sum\limits_{g\in\mathcal{G}_{r(k)}}\phi_{i}^{p}(r(k))\left|\mu_{i}(g)-\mu_{i}(gk)\right|+\mu_{i}(gk)\left|\phi_{i}^{p}(r(k))-\phi_{i}^{p}(s(k))\right|\\
&=\sum\limits_{i=0}^{f([(R+1)^{c}]R+1)}\phi_{i}^{p}(r(k))\cdot\sum\limits_{g\in(\mathcal{G}_{i})_{r(k)}}\left|\mu_{i}(g)-\mu_{i}(gk)\right|\\&+\sum\limits_{i=0}^{f([(R+1)^{c}]R+1)}\left|\phi_{i}(r(k))-\phi_{i}(s(k))\right|\cdot\left|\sum\limits_{j=0}^{p-1}\phi_{i}(r(k))^{j}\phi_{i}(s(k))^{p-1-j}\right|\cdot\sum\limits_{g\in(\mathcal{G}_{i})_{r(k)}}\mu_{i}(gk)\\&<\sum\limits_{i=0}^{f([(R+1)^{c}]R+1)}\phi_{i}^{p}(r(k))\cdot\frac{\varepsilon}{2}+p\cdot\sum\limits_{i=0}^{f([(R+1)^{c}]R+1)}\left|\phi_{i}(r(k))-\phi_{i}(s(k))\right|\\
&<\frac{\varepsilon}{2}+\frac{\varepsilon}{2}=\varepsilon, 
\end{align*}
where the first equality uses the equality $\left|a^{p}-b^{p}\right|=\left|a-b\right|\cdot\left|\sum\limits_{j=0}^{p-1}a^{j}b^{p-1-j}\right|$ for all $a, b\geq 0$, the third inequality is due to the amenability of $\mathcal{G}_{i}$, and the final inequality follows from the item (3) of Proposition \ref{prop 5.1}.
\end{proof}

\begin{remark}
In the final steps of the third part of the proof, to obtain the inequality 
\begin{equation*}
\sum\limits_{g\in\mathcal{G}_{r(k)}}\left|\mu(g)-\mu(gk)\right|<\varepsilon,
\end{equation*}
we need to show that 
\begin{align}\label{eq 6.1}
p\cdot\sum\limits_{i=0}^{f([(R+1)^{c}]R+1)}\left|\phi_{i}(r(k))-\phi_{i}(s(k))\right|<\frac{\varepsilon}{2}.
\end{align}
As shown in the final steps of the proof of Proposition \ref{prop 5.1}, the left-hand side of inequality \eqref{eq 6.1} is bounded above by the left-hand side of inequality \eqref{eq 5.1}. Hence, \eqref{eq 6.1} holds provided that the left-hand side of \eqref{eq 5.1} tends to zero as $R\to\infty$. This condition implies that $$\frac{f(N(R)R)^{\frac{1}{p}+1}}{N(R)}\to 0, \text{ as }R\to\infty.$$ Consequently, the growth rate of $f$ must be slower than linear; for example, $f(x)=[x^{\alpha}]$ with $0<\alpha<1$ satisfies this requirement. In this case, choosing $p=[\frac{\alpha}{1-\alpha}]+1$ and $N(R)=[(R+1)^{c}]$ for a suitable constant $c$ (see Remark \ref{rmk 5.7}) ensures that $f(N(R)R)^{\frac{1}{p}+1}$ grows slower than $N(R)$. However, if $f$ satisfies stronger growth conditions, such as linear or superlinear growth, the above limit is no longer zero, and the required inequality cannot be obtained by this approach. Therefore, the current proof of Theorem \ref{th 6.4} fails to hold under stronger growth conditions.
\end{remark}

The following is an example demonstrating that slow dynamic asymptotic dimension growth of a groupoid implies its amenability, obtained through a combination of Example \ref{example 2.11}, Theorem \ref{th 4.17}, and Theorem \ref{th 6.4}.
\begin{example}\label{example 6.6}
Let $f: \mathbb{N}\rightarrow\mathbb{R}^{+}$ be a strictly increasing function such that $f\succeq x^{\alpha}$ for some $\alpha>1$. Consider the space $X=\bigoplus\limits_{n=1}^{\infty}f(n)\mathbb{Z}$, equipped with the direct sum metric given in \eqref{eq 2.2}. Then $X$ is a discrete metric space with bounded geometry and slow asymptotic dimension growth $f^{-1}$. According to Theorem \ref{th 4.17}, the corresponding coarse groupoid $G(X)$ has slow dynamic asymptotic dimension growth $f^{-1}\preceq x^{\frac{1}{\alpha}}$ with $0<\frac{1}{\alpha}<1$. Consequently, Theorem \ref{th 6.4} implies that $G(X)$ is amenable. 
\end{example}

This section concludes by investigating the implications of dynamic asymptotic dimension growth for two major conjectures: the Baum--Connes conjecture and Matui's HK conjecture. Finally, we propose several open problems arising from our study.

Let $\mathcal{G}$ be a $\sigma$-compact, locally compact, Hausdorff groupoid with Haar system, and let $B$ be a $\mathcal{G}$-algebra. The Baum–Connes assembly map is defined as the homomorphism
$$
\mu_{r}^{B}: K^{\mathrm{top}}_{*}(\mathcal{G};B)\rightarrow K_{*}(B\rtimes_{r} \mathcal{G}).
$$ 
The Baum--Connes conjecture with coefficients in $B$ asserts that $\mu_{r}^{B}$ is an isomorphism.

By combining Theorem \ref{th 6.4} with the established result from \cite{Tu99} that the Baum--Connes conjecture with coefficients holds for amenable groupoids, we derive the following corollary:

\begin{corollary}\label{cor 6.6}
Let $\mathcal{G}$ be a $\sigma$-compact, locally compact, Hausdorff, and \'etale groupoid with compact unit space and dynamic asymptotic dimension growth $f\preceq x^{\alpha} (0<\alpha<1)$. Then the Baum--Connes conjecture with coefficients holds for $\mathcal{G}$.
\end{corollary}

The following example follows directly from Example \ref{example 6.6} and Corollary \ref{cor 6.6}. 
\begin{example}\label{example 6.8}
Let $f: \mathbb{N}\rightarrow\mathbb{R}^{+}$ be a strictly increasing function such that $f\succeq x^{\alpha}$ for some $\alpha>1$, and let $X=\bigoplus\limits_{n=1}^{\infty}f(n)\mathbb{Z}$ be a space as in Example \ref{example 6.6}. Then the Baum--Connes conjecture with coefficients holds for the coarse groupoid $G(X)$.
\end{example}

The following theorem, [\cite{PY26}, Corollary 5.6], demonstrates that the rational version of Matui's HK conjecture holds for a large class of ample \'etale groupoids.

\begin{theorem}[\cite{PY26}]\label{th 6.9}
Let $\mathcal{G}$ be a second-countable, locally compact, Hausdorff, ample, and \'etale groupoid with torsion-free isotropy. If $\mathcal{G}$ satisfies the rational Baum--Connes conjecture, then the rational HK conjecture holds for $\mathcal{G}$, i.e., we have the following isomorphism of abelian groups:
\begin{equation*}
K_{i}(C^{*}_{r}(\mathcal{G}))\otimes\mathbb{Q}\cong\bigoplus_{k\in\N}H_{i+2k}(\mathcal{G};\mathbb{Q}) \quad(i=0, 1).
\end{equation*}
\end{theorem}

By combining Corollary \ref{cor 6.6} with Theorem \ref{th 6.9}, we obtain the following corollary, which establishes a sufficient condition for the rational HK conjecture. 
\begin{corollary}\label{cor 6.10}
Let $\mathcal{G}$ be a second-countable, $\sigma$-compact, locally compact, Hausdorff, ample, and \'etale groupoid with compact unit space and torsion-free isotropy. If $\mathcal{G}$ has dynamic asymptotic dimension growth $f\preceq x^{\alpha} (0<\alpha<1)$, then the rational HK conjecture holds for $\mathcal{G}$.
\end{corollary}

The following example is a consequence of Example \ref{example 6.8}, [\cite{STY02}, Section 3.3] and Corollary \ref{cor 6.10}. 
\begin{example}
Let $G(X)$ be the coarse groupoid as in Example \ref{example 6.8}. According to [\cite{STY02}, Section 3.3], there exists a subgroupoid $\mathcal{G}$ of $G(X)$ such that $\mathcal{G}$ is second-countable, $\sigma$-compact, locally compact, Hausdorff, ample, principal, and \'etale, with compact unit space and dynamic asymptotic dimension growth at most $x^{\alpha}(0<\alpha<1)$. Hence, the rational HK conjecture holds for $\mathcal{G}$.
\end{example}

Having established these results, we now formulate several open questions. A natural question is the following, motivated by Tu's result that the Baum--Connes conjecture with coefficients holds for a-T-menable groupoids (i.e., those satisfying the Haagerup property) \cite{Tu99}:

\begin{question}\label{q 6.12}
What is the relationship between the dynamic asymptotic dimension growth of a groupoid and its a-T-menability?
\end{question}

It is established that a groupoid $\mathcal{G}$ has the Haagerup property (i.e., is a-T-menable) if it admits a continuous proper negative type function \cite{Tu99}. Therefore, answering Question \ref{q 6.12} depends on establishing a connection between dynamic asymptotic dimension growth and the existence of such a function, for which methods inspired by \cite{Tu99, STY02} may be applicable.

A positive answer to Question \ref{q 6.12}, combined with the result that a-T-menable groupoids satisfy the Baum--Connes conjecture \cite{Tu99}, would provide a pathway to address the central question for \'etale groupoids:

\begin{question}\label{q 6.13}
What kind of dynamic asymptotic dimension growth of \'etale groupoids implies the Baum--Connes conjecture?
\end{question}

In particular, we expect that polynomial dynamic asymptotic dimension growth implies the Baum--Connes conjecture. This expectation is motivated by the following chain of implications: slow dynamic asymptotic dimension growth implies amenability, which in turn implies the conjecture through the a-T-menability condition---a property weaker than amenability. Therefore, we anticipate that a stronger growth condition, such as polynomial growth, may also suffice to establish the conjecture.

Finally, in the context of second-countable, principal, ample, and \'etale groupoids, we pose the following question:

\begin{question}
For a second-countable, locally compact, Hausdorff, principal, ample, and \'etale groupoid, what kind of dynamic asymptotic dimension growth implies the rational HK conjecture?
\end{question}

The answer to Question \ref{q 6.13}, together with results such as Theorem \ref{th 6.9}, may provide insights into the final question.

\section*{Acknowledgments}
We thank Jiawen Zhang, Kang Li, Zehong Huang, Xin Ma, Jianchao Wu, and Valerio Proietti for their valuable discussions and suggestions.

\end{document}